\magnification=1202
\input amssym1.def
\font\bb=msbm10 at9.98pt
\font\bbm=bbm10 at9.98pt
 at9.98pt
 at9.98pt
\font\callig=callig15 at15pt
 at9.98pt 
\font\cyr=wncyi10 at 9.98pt
\font\eightrm=cmr8
\font\eightsc=cmcsc8

 at9.98pt
\font\larm=cmr10 at13.1pt
\font\labf=cmbx10 at13.1pt
 at15.74pt
\font\rsfs=rsfs10 at9.98pt
\font\scriptrsfs=rsfs7 at6.98pt
\font\scriptscriptrsfs=rsfs5 at 5pt
 at7.6pt
 at9.98pt
 at9.98pt
 at9.98pt
 at9.98pt
\font\sc=cmcsc10 at9.98pt
 at8.31pt
\font\tenpbf=cmbx10 at8.32pt
\font\tenpit=cmti10 at8.32pt
\font\tenprm=cmr10 at8.32pt
 at9.1pt
 at10.92pt
\textfont4=\rsfs \scriptfont4=\scriptrsfs 
\scriptscriptfont4=\scriptscriptrsfs
\def\script{\fam4 \rsfs}


\def\meriki{\hbox{\cyr\char'144}}
\def\Der{\hbox{$\yfra D\yfra e\yfra r$}}

\def\sT*MA{\hbox{${\script T}_{p}^{\ast}(M,{\script A})$}}
\def\kenosyn{\hbox{\bb\char'77}}


\def\A{\hbox{$\script A$}}
\def\Adneo{{\rm Ad\kern1pt}}

\def\an{\raise0.5pt\hbox{$\kern2pt\scriptstyle\in\kern2pt$}}
\def\arkef{\advance\chapternumber by 1\sc\rm{\the\chapternumber}}
\def\Auto{{\rm Aut\kern1pt}}
\def\B{\hbox{$\script B$}}

\def\C{\hbox{$\script C$}}
\def\Coad{{\rm Coad\kern1pt}}

\def\Colon{\colon\kern2pt}
\def\comp{\hbox{\lower5.8pt\hbox{\larm\char'027}}}
\def\crossed{\hbox{\bb\char'076}}
\def\D{\hbox{$\script D$}}
\def\da{\hbox{$\Delta_{\script A}$}}
\def\dacirc{\hbox{$\Delta^{\circ}_{\script A}$}}

\def\dim{{\rm dim\kern1pt}}
\def\deg{{\rm deg\kern1pt}}
\def\E{\hbox{$\script E$}}
\def\ena{\hbox{\bbm 1}}
\def\Endo{{\rm End\kern1pt}}
\def\eightpoint{\eightrm}
\def\F{\hbox{$\script F$\kern1pt}}
\catcode`\@=11
\def\footnote#1{\edef\@sf{\spacefactor\the\spacefactor}#1\@sf
     \insert\footins\bgroup\eightpoint
     \interlinepenalty100 \let\par=\endgraf
      \leftskip=0pt \rightskip=0pt
      \splittopskip=10pt plus 1pt minus 1pt \floatingpenalty=20000
      \smallskip\item{#1}\bgroup\strut\aftergroup\@foot\let\neft}
\skip\footins=12pt plus 2pt minus 4pt
\dimen\footins=30pc
\def\H{\hbox{$\script H$}}
\def\hom{\hbox{\callig h\kern5.5pt}}
\def\im{{\rm im\kern1pt}}

\def\ker{{\rm ker\kern1pt}}
\def\Lie{\hbox{$\bit L$}}
\def\line{\hbox to\hsize}
\def\meriki{\hbox{\cyr\char'144\kern0.3pt}}

\def\na{\raise0.5pt\hbox{$\kern2pt\scriptstyle\ni\kern2pt$}}
\def\noan{\hbox{$\an\raise0.5pt\hbox{$\kern-6.5pt\scriptstyle\slash\kern3pt$}$}}
\def\O{\hbox{$\script O$}}
\def\odot{\;{\mathchar"220C}\;}

\def\opi{\hbox{$\kern-1pt\hat{\otimes}_{\hskip-0.1cm\pi}\kern1pt$}}
\def\oplus{\;{\mathchar"2208}\;}

\def\otimes{\;{\mathchar"220A}\;}
\def\OYBG{\hbox{$\Omega(Y,\B,{\frak g})$}}
\def\R{\hbox{\bf\char'122}}
\def\Pro{\hbox{$\script P$\kern-2pt\callig ro\kern2pt}}
\def\pounds{\rlap{\lower3.5pt\hbox{\kern2.9pt\hbox{\char'26}}}
{\script L}}
\def\qed{\hbox{\kern0.3cm\vrule height6.7pt width6.7pt depth-0.2pt}}
\def\qedaux{\hbox{\kern0.4cm\vrule height4.5pt width4.5pt
depth-0.2pt}}
\def\san{\raise0.5pt\hbox{$\kern0.7pt\scriptscriptstyle\in\kern0.7pt$}}
\def\scomp{\hskip-0.05truecm\hbox{\lower5pt\hbox{$\mathchar"2017$}}
\hskip-0.05truecm}

\def\slje{\hbox{$\meriki\kern-4.6pt\raise0.7pt\hbox{\slash}$}}
\def\swH{\kern4pt\hbox{\kern4pt\hbox{$\widetilde{
}$}\kern-7pt{$\scriptstyle\script H$}}}
\def\T{\hbox{$\script T$\kern1.5pt}}
\def\times{\;{\mathchar"2202}\;}
\def\tonos{\hbox{\kern-1.3pt\lower0.7pt\hbox{$\mathchar"6013$}}}
\def\tonoskef{\hbox{\kern-1.3pt\hbox{$\mathchar"6013$}}}
\def\tonoss{\hbox{\kern-1.3pt\lower0.9pt\hbox{\tenprm\char'23}}}
\def\Ver{\hbox{$\script V$\kern-3.5pt\callig er\kern2pt}}
\def\wE{\hbox{$\widetilde{\E\ }\kern-1pt$}}
\def\wi{\hbox{$\widetilde{\,\,\,}$}}
\def\wH{\kern2.5pt\hbox{\raise2.5pt\hbox{\wi}\kern-16pt{\H}}}

\def\,{\hskip0.2truecm}
\def\Z{\hbox{\bf\char'132}}

                                               %
\def\author#1{{\tenprm #1:}}                   %
\def\ekdoths#1{{\tenprm #1}}                   %
\def\periodiko#1{{\tenpit #1\tenprm ,}}	       %
\def\titlosa#1{{\tenprm #1,}}                  %
\def\titlosb#1{{\tenpit #1\tenprm ,}}          %
\def\selides#1{{\tenprm #1}}		       %
\def\tomos#1{{\tenprm Vol. \tenpbf #1\tenprm :}}%
					       %

\def\title#1{\line{\hss}\line{\hss#1\hss}%
\line{\hss}\hskip-0.75truecm}

%
%
%
\def\teleia{\hbox{.}}
\newif\ifPhysRev
\def\Textindent#1{\noindent\llap{#1\enspace}\ignorespaces}
\def\nonfrenchspacing{\sfcode`\.=3001 \sfcode`\!=3000 \sfcode`\?=3000
        \sfcode`\:=2000 \sfcode`\;=1500 \sfcode`\,=1251 }
\nonfrenchspacing
\newdimen\d@twidth
 {\setbox0=\hbox{s.} \global\d@twidth=\wd0 \setbox0=\hbox{s}
        \global\advance\d@twidth by -\wd0 }
\def\removehglue{\loop \unskip \ifdim\lastskip >\z@ \repeat }
\def\roll@ver#1{\removehglue \nobreak \count255 =\spacefactor \dimen@=\z@
        \ifnum\count255 =3001 \dimen@=\d@twidth \fi
        \ifnum\count255 =1251 \dimen@=\d@twidth \fi
    \iftwelv@ \kern-\dimen@ \else \kern-0.83\dimen@ \fi
   #1\spacefactor=\count255 }
\def\step@ver#1{\relax \ifmmode #1\else \ifhmode
        \roll@ver{${}#1$}\else {\setbox0=\hbox{${}#1$}}\fi\fi }
\def\attach#1{\step@ver{\strut^{\mkern 2mu #1} }}

\normalbaselineskip = 20pt plus 0.2pt minus 0.1pt
\normallineskip = 1.5pt plus 0.1pt minus 0.1pt
\normallineskiplimit = 1.5pt
\newskip\normaldisplayskip
\normaldisplayskip = 20pt plus 5pt minus 10pt
\newskip\normaldispshortskip
\normaldispshortskip = 6pt plus 5pt
\newskip\normalparskip
\normalparskip = 6pt plus 2pt minus 1pt
\newskip\skipregister
\skipregister = 5pt plus 2pt minus 1.5pt
\newif\ifsingl@    \newif\ifdoubl@
\newif\iftwelv@    \twelv@true
\def\singlespace{\singl@true\doubl@false\spaces@t}
\def\doublespace{\singl@false\doubl@true\spaces@t}
\def\normalspace{\singl@false\doubl@false\spaces@t}
\def\Tenpoint{\tenpoint\twelv@false\spaces@t}
\def\Twelvepoint{\twelvepoint\twelv@true\spaces@t}
\def\spaces@t{\relax
      \iftwelv@ \ifsingl@\subspaces@t3:4;\else\subspaces@t1:1;\fi
       \else \ifsingl@\subspaces@t3:5;\else\subspaces@t4:5;\fi \fi
      \ifdoubl@ \multiply\baselineskip by 5
         \divide\baselineskip by 4 \fi }
\def\subspaces@t#1:#2;{
      \baselineskip = \normalbaselineskip
      \multiply\baselineskip by #1 \divide\baselineskip by #2
      \lineskip = \normallineskip
      \multiply\lineskip by #1 \divide\lineskip by #2
      \lineskiplimit = \normallineskiplimit
      \multiply\lineskiplimit by #1 \divide\lineskiplimit by #2
      \parskip = \normalparskip
      \multiply\parskip by #1 \divide\parskip by #2
      \abovedisplayskip = \normaldisplayskip
      \multiply\abovedisplayskip by #1 \divide\abovedisplayskip by #2
      \belowdisplayskip = \abovedisplayskip
      \abovedisplayshortskip = \normaldispshortskip
      \multiply\abovedisplayshortskip by #1
        \divide\abovedisplayshortskip by #2
      \belowdisplayshortskip = \abovedisplayshortskip
      \advance\belowdisplayshortskip by \belowdisplayskip
      \divide\belowdisplayshortskip by 2
      \smallskipamount = \skipregister
      \multiply\smallskipamount by #1 \divide\smallskipamount by #2
      \medskipamount = \smallskipamount \multiply\medskipamount by 2
      \bigskipamount = \smallskipamount \multiply\bigskipamount by 4 }
\def\normalbaselines{ \baselineskip=\normalbaselineskip
   \lineskip=\normallineskip \lineskiplimit=\normallineskip
   \iftwelv@\else \multiply\baselineskip by 4 \divide\baselineskip by 5
     \multiply\lineskiplimit by 4 \divide\lineskiplimit by 5
     \multiply\lineskip by 4 \divide\lineskip by 5 \fi }


\def\abstract#1{\parshape=1 0.7cm \dimen10
                {\tenpbf Abstract. \tenprm #1}}

\newcount\appendixnumber     \appendixnumber=0
\newcount\chapternumber      \chapternumber=0
\newcount\equanumber         \equanumber=0
\newcount\mathnumber         \mathnumber=0
\newcount\appequanumber      \appequanumber=0
\newcount\appmathnumber      \appmathnumber=0

\let\variableone=\relax
\let\variabletwo=\relax
\let\chapterlabel=\relax
\let\sectionlabel=\relax
\let\mathlabel=\relax
\newtoks\chapterstyle        \chapterstyle={\Number}
\newtoks\sectionstyle        \sectionstyle={\chapterlabel\Number}
\newskip\chapterskip         \chapterskip=\bigskipamount
\newskip\sectionskip         \sectionskip=\medskipamount
\newskip\headskip            \headskip=8pt plus 3pt minus 3pt
\newdimen\chapterminspace    \chapterminspace=15pc
\newdimen\sectionminspace    \sectionminspace=10pc
\newdimen\sectionspace       \sectionspace=20pc
\newdimen\referenceminspace  \referenceminspace=25pc

\def\chapterreset{\global\advance\chapternumber by 1
   \ifnum\equanumber<0 \else\global\equanumber=0\fi
   \mathnumber=0
   \makechapterlabel}
\def\makechapterlabel{\let\sectionlabel=\relax\let\mathlabel=\relax
 \xdef\chapterlabel{\the\chapterstyle{\the\chapternumber\teleia\kern3pt}}}

\def\rightheadline{\sc\hfil\variableone\eightsc\hfil\folio}
\def\leftheadline{\eightsc\folio\hfil{\sc\variabletwo}\hfil}
\def\heads{\footline={\hfil}\headline={\ifodd\pageno
               \rightheadline\else\leftheadline\fi}}

\def\headseis{\partreset\headline={\ifodd\pageno{
                         \hfil\sc partie {\eightsc\the\partnumber}
                         -introduction\hfil\eightsc\folio}\else
                        {\eightsc\folio\hfil\sc partie
                         {\eightsc\the\partnumber}-introduction\hfil}\fi}
                        \footline={\hfil}}

\def\alphabetic#1{\count255='140 \advance\count255 by #1\char\count255}
\def\Alphabetic#1{\count255='100 \advance\count255 by #1\char\count255}
\def\Roman#1{\uppercase\expandafter{\romannumeral #1}}

\def\Number#1{\number #1}
\def\BLANC#1{}

\def\titlestyle#1{\par\begingroup \interlinepenalty=9999
     \leftskip=0.02\hsize plus 0.23\hsize minus 0.02\hsize
     \rightskip=\leftskip \parfillskip=0pt
     \hyphenpenalty=9000 \exhyphenpenalty=9000
     \tolerance=9999 \pretolerance=9000
     \spaceskip=0.333em \xspaceskip=0.5em
     \iftwelv@\bf\else\bf\fi
   \noindent #1\par\endgroup }

\def\spacecheck#1{\dimen@=\pagegoal\advance\dimen@ by -\pagetotal
   \ifdim\dimen@<#1 \ifdim\dimen@>0pt \vfil\break \fi\fi}
\def\TableOfContentEntry#1#2#3{\relax}

\def\chapter#1{\vskip0.5cm

   \par\penalty -300 \vskip\chapterskip
   \chapterreset \titlestyle{\chapterlabel\ #1}
   \nobreak\vskip\headskip \penalty 30000
   \wlog{\string\chapter\space \chapterlabel} }

\def\appendixreset{\global\advance\appendixnumber by 1
                   \appmathnumber=0\appequanumber=0}
\def\appendix#1{\par \penalty-300\vskip\chapterskip
   \spacecheck\chapterminspace
   \appendixreset \title{\bf Appendix \Alphabetic{\the\appendixnumber}}
   \nobreak\vskip-\chapterskip\penalty 30000
   \vskip-\chapterskip
   \par{\titlestyle{#1}}
   \vskip\chapterskip
   \wlog{\string\appendix\space \chapterlabel} }

%
%
\def\eqname#1{\relax \ifnum\equanumber<0
     \xdef#1{{\noexpand\rm(\number-\equanumber)}}%
       \global\advance\equanumber by -1
    \else \global\advance\equanumber by 1
      \xdef#1{{\noexpand(\rm{\the\chapternumber}\teleia
                            \rm{\number\equanumber})}} \fi #1}

\def\eqn{\eqno\eqname}

\def\math#1#2{\vskip0.1cm
   \global\advance\mathnumber by 1
   \xdef\mathlabel{\the\chapternumber\teleia\the\mathnumber}
   \wlog{\string\math\space \mathlabel}
   {\bf\enspace\mathlabel\hskip0.2cm #1}
   \xdef#2{{\mathlabel}}}

\def\appeqname#1{\relax \ifnum\appequanumber<0
     \xdef#1{{\noexpand\rm(\number-\appequanumber)}}%
       \global\advance\appequanumber by -1
    \else \global\advance\appequanumber by 1
      \xdef#1{{\noexpand(\hbox{\Alphabetic{\the\appendixnumber}}\teleia
                            {\number\appequanumber})}} \fi #1}

\def\mathapp#1#2{\vskip0.1cm
   \global\advance\appmathnumber by 1
   \xdef\appmathlabel{{\Alphabetic{\the\appendixnumber}}\teleia
   \the\appmathnumber}
   \wlog{\string\mathapp\space \appmathlabel}
   {\bf\enspace\appmathlabel\hskip0.2cm #1}
   \xdef#2{{\appmathlabel}}}


%
%
%
\newtoks\referencestyle      \referencestyle={\tenpbf\Number}
\newcount\referencecount     \referencecount=0
\newcount\lastrefsbegincount \lastrefsbegincount=0
\newif\ifreferenceopen       \newwrite\referencewrite
\newif\ifrw@trailer
\newdimen\refindent     \refindent=12pt
\def\NPrefmark#1{\attach{\scriptscriptstyle [ #1 ] }}
\let\PRrefmark=\attach
\def\refmark#1{\relax\ifPhysRev\PRrefmark{#1}\else\NPrefmark{#1}\fi}
\def\refend@{\refmark{\number\referencecount}}
\def\refend{\refend@{}\space }
\def\refsend{\refmark{\count255=\referencecount
   \advance\count255 by-\lastrefsbegincount
   \ifcase\count255 \number\referencecount
   \or \number\lastrefsbegincount,\number\referencecount
   \else \number\lastrefsbegincount-\number\referencecount \fi}\space }
\def\refitem#1{\par\hangafter=0 \hangindent=\refindent	\Textindent{#1}}
\def\Ref{\rw@trailertrue\REF}
\def\REF#1{\r@fstart{#1}%
   \rw@begin{\tenprm [\tenpbf\Number{\the\referencecount}\tenprm ]}\rw@end}
\def\r@fstart#1{\chardef\rw@write=\referencewrite \let\rw@ending=\refend@
   \ifreferenceopen \else \global\referenceopentrue
   \immediate\openout\referencewrite=referenc.txa
   \toks0={\catcode`\^^M=10}\immediate\write\rw@write{\the\toks0} \fi
   \global\advance\referencecount by 1 
   \xdef#1{[{\the\referencestyle{\the\referencecount}}]}}
 {\catcode`\^^M=\active %
 \gdef\rw@begin#1{\immediate\write\rw@write{\noexpand\refitem{#1}}%
   \begingroup \catcode`\^^M=\active \let^^M=\relax}%
 \gdef\rw@end#1{\rw@@end #1^^M\rw@terminate \endgroup%
   \ifrw@trailer\rw@ending\global\rw@trailerfalse\fi }%
 \gdef\rw@@end#1^^M{\toks0={#1}\immediate\write\rw@write{\the\toks0}%
   \futurelet\n@xt\rw@test}%
 \gdef\rw@test{\ifx\n@xt\rw@terminate \let\n@xt=\relax%
       \else \let\n@xt=\rw@@end \fi \n@xt}%
}
\let\rw@ending=\relax
\let\rw@terminate=\relax

\def\Closeout#1{\toks0={\catcode`\^^M=5}\immediate\write#1{\the\toks0}%
   \immediate\closeout#1}

\input pictex 
\topskip1truecm
\voffset=2.5truecm
\hsize 15truecm
\vsize 20truecm
\hoffset=0.25truecm
\def\undertext#1{$\underline{\hbox{#1}}$}
\topglue 3cm
\pageno=1
\baselineskip=14pt plus .2pt
\footline={\hss\eightsc\folio\hss}
\dimen10=\hsize \advance\dimen10 by -1.4cm
\def\variableone{t. stavracou}
\def\variabletwo{theory of connections on graded principal bundles}


\Ref\alm{
\author{A. L. Almorox}
\titlosa{Supergauge theories in graded manifolds}
\periodiko{Differential Geometrical Methods in Mathematical Physics}
\ekdoths{Proceedings, Salamanca 1985, Eds.:
P.L. Gar\-c\'{\i}a, A. P\'erez-Rend\'on, Springer LNM 1251}}

\Ref\almdyo{
\author{A. L. Almorox}
\titlosa{The bundle of graded frames}
\periodiko{Rend. Sem. Mat. Univers. Politech. Torino}
\tomos{43}
\selides{405--426 (1985)}}

\Ref\bbr{
\author{C. Bartocci, U. Bruzzo, D. H. Ruip\'erez}
\titlosb{The geometry of supermanifolds}
\ekdoths{Kluwer Academic Publishers, 1991}}

\Ref\batchelora{
\author{M. Batchelor}
\titlosa{The structure of supermanifolds}
\periodiko{Trans. Amer. Math. Soc.}
\tomos{253}
\selides{329--338 (1979),}
\titlosa{Two approaches to supermanifolds} 
\periodiko{Trans. Amer. Math. Soc.}
\tomos{258}
\selides{257--270 (1980)}}

\Ref\berezin{
\author{F. Berezin, D. Leites}
\titlosa{Supermanifolds} 
\periodiko{Soviet Math. Dokl.}
\tomos{16}
\selides{1218--1222 (1975)}}

\Ref\bmarinov{
\author{F. Berezin, M. Marinov}
\titlosa{Particle spin dynamics as the Grassmann variant of 
classical mechanics} 
\periodiko{Ann. of Phys.}
\tomos{104}
\selides{336--362 (1977)}}

\Ref\bruzzo{
\author{U. Bruzzo, R. Cianci}
\titlosa{Mathematical theory of super fibre bundles}
\periodiko{Class. Quantum Grav.}
\tomos{1}
\selides{213--226 (1984)}}

\Ref\cant{
\author{F. Cantrijn, L.A. Ibort}
\titlosa{Introduction to Poisson supermanifolds}
\periodiko{Diff. Geom. Appl.}
\tomos{1}
\selides{133--152 (1991)}}

\Ref\riccigia{
\author{R. Giachetti, R. Ricci}
\titlosa{{\tenpbf R}-actions, derivations, and Frobenius theorem on 
graded manifolds}
\periodiko{Adv. Math.}
\tomos{62}
\selides{84--100 (1986)}}

\Ref\grasso{
\author{M. Grasso, P. Teofilatto}
\titlosa{Gauge theories, flat superforms and reduction of super fibre
bundles}
\periodiko{Rep. Math. Phys.}
\tomos{25}
\selides{53--71 (1987)}}

\Ref\grot{
\author{A. Grothendieck}
\titlosa{Produits tensoriels topologiques et espaces nucl\'eaires}
\periodiko{Mem. A.M.S.}
\tomos{16}
\selides{(1955)}}

\Ref\ibortsolano{
\author{A. Ibort, J. Mar\'{\i}n-Solano}
\titlosa{Geometrical foundations of Lagrangian supermechanics and 
supersymmetry}
\periodiko{Rep. Math. Phys.}
\tomos{32}
\selides{385--409 (1993)}}

\Ref\arcad{
\author{A. Jadczyk, K. Pilch}
\titlosa{Superspaces and supersymmetries}
\periodiko{Comm. Math. Phys.}
\tomos{78}
\selides{373--390 (1981)}}

\Ref\kobay{
\author{S. Kobayashi, K. Nomizu}
\titlosb{Foundations of differential geometry}
\ekdoths{Interscience tracts in pure and applied mathematics, John
Wiley and Sons, 1963}}

\Ref\kost{
\author{B. Kostant}
\titlosa{Graded manifolds, graded Lie theory and 
prequantization}
\periodiko{Differential geometric methods in 
mathematical physics}
\ekdoths{Lecture Notes in Mathematics}
\tomos{570}
\selides{177--306,}
\ekdoths{Berlin, Heidelberg, New York: Springer-Verlag, 1977}}

\Ref\leites{
\author{D. Leites}
\titlosa{Introduction to the theory of supermanifolds}
\periodiko{Russ. Math. Surv.}
\tomos{35}
\selides{1--64 (1980)}}

\Ref\montedyo{
\author{J. Monterde}
\titlosa{Higher order graded and Berezinian Lagrangian densities 
and their Euler-Lagran\-ge equations}
\periodiko{Ann. Inst. Henri Poincar\'e (Phys. Th\'eor.)}
\tomos{57}
\selides{3--26 (1992)}}

\Ref\monterma{
\author{J. Monterde, J. M. Masqu\'e}
\titlosa{Variational problems on graded manifolds}
\periodiko{Contemp. Math.}
\tomos{132}
\selides{551--571 (1992)}}

\Ref\mont{
\author{S. Montgomery}
\titlosb{Hopf algebras and their actions on rings}
\ekdoths{CBMS, Regional Conference Series in Mathematics, Number 82,
American Mathematical Society, 1993}}

\Ref\rrr{
\author{A. P\'erez-Rend\'on, D.H. Ruip\'erez}
\titlosa{Towards a classical field theory on graded manifolds}
\periodiko{Proc. ``Journ\'ees relativistes"}
\ekdoths{Torino, 1983}}

\Ref\roger{
\author{A. Rogers}
\titlosa{A global theory of supermanifolds}
\periodiko{J. Math.  Phys.}
\tomos{21}
\selides{1352--1365 (1980)}}
 
\Ref\rogers{
\author{A. Rogers}
\titlosa{Graded manifolds, supermanifolds and infinite-dimensional
Grassmann algebras}\break
\periodiko{Comm. Math. Phys.}
\tomos{105}
\selides{375--384 (1986)}}

\Ref\ruipma{
\author{D. H. Ruip\'erez, J. Mu\~noz Masqu\'e}
\titlosa{Global variational calculus on graded manifolds I: graded jet
bundles, structure 1-forms and graded infinitesimal contact
transformations}
\periodiko{J. Math. Pures et Appl.}
\tomos{63}
\selides{283--309 (1984) and}
\titlosa{Global variational calculus on graded manifolds II}
\periodiko{J. Math. Pures et Appl.}
\tomos{63}
\selides{87--104 (1985)}}

\Ref\dewitt{ 
\author{B. DeWitt}
\titlosb{Supermanifolds}
\ekdoths{Cambridge University Press, 1992}}

\vskip-2.5cm
\centerline{\bf Centre de Physique Th\'eorique\footnote{$^\ast$}{Unit\'e 
Propre de Recherche 7061} - CNRS - Luminy, Case 907}
\centerline{\bf F-13288 Marseille Cedex 9 - France}
\vskip1.5cm
\centerline{\labf Theory of connections on graded principal
bundles}
\vskip0.7cm
\centerline{\bf T. Stavracou\footnote{$^{1}$}{email: 
jenny@cpt.univ-mrs.fr}}
\vskip1cm

\advance\baselineskip by -1pt

\parshape=1 0.7cm \dimen10{
{{\bf Abstract.} The geometry of graded principal bundles is discussed 
in the framework of graded manifold theory of Kostant-Berezin-Leites. 
In particular, we prove that a graded principal bundle is globally 
trivial if and only if it admits a global graded section and, further,
that the sheaf of vertical derivations on such a bundle coincides with 
the graded distribution induced by the action of the structure 
graded Lie group. This result leads to a natural definition of the 
graded connection in terms of graded distributions; its relation 
with Lie superalgebra-valued graded differential forms is also exhibited. 
Finally, we define the curvature for the graded connection; we prove 
that the curvature controls the involutivity of the horizontal graded 
distribution corresponding to the graded connection.}}

\advance\baselineskip by 1pt

\vskip0.8cm
\parshape=1 0.7cm \dimen10{

\vskip0.1cm


\vskip0.1cm

April 1996

CPT-96/P.3331

\vskip0.2cm

\noindent anonymous ftp or gopher: cpt.univ-mrs.fr}

\vfill\eject

\chapter{Introduction}
\heads

Graded manifolds are geometrical objects introduced by Kostant, 
{\kost} and Berezin, Leites, {\berezin}, as the natural 
mathematical tool in order to study supersymmetric problems. In 
particular, the formalism of graded manifolds provides the possibility
to deal with classical dynamics of particles with spin. Indeed, it has
been proved, {\bmarinov}, that the extension of ordinary phase spaces
by anti-commuting variables yields the classical phase spaces of
particles with spin. Thus, Kostant's article, {\kost}, can be also
considered as a treatment on super-Hamiltonian systems and on their
geometrical prequantization. We mention {\ibortsolano}, {\montedyo}, 
{\monterma}, {\ruipma} as most important references related more 
or less to the global problems of Lagrangian supermechanics. 

Our aim here is to investigate several aspects of the graded
differential geometry, unexplored up today, and to use them in order
to establish a theory of connections on graded principal bundles. 
This work is mainly motivated by the need to make clear why graded
connection theory is convenient for the description of supersymmetric 
gauge theories, and how one could use this theory in practice, in
order to obtain, for example, the supersymmetric standard model
in the same spirit as one obtains the standard model of fundamental
interactions from ordinary gauge theory. The notion of graded 
principal bundle has first appeared in {\alm}, {\almdyo}: there, 
the notion of connection on such bundles has been also presented, 
but in a rather academic way, i.e. using jet bundles, bundles of 
connections and Atiyah's exact sequence, conveniently generalized 
in this context. This point of view turned out to be particularly 
useful for mathematics, as it inspired and affected more recent 
developments in connection theory on far more sophisticated types 
of supermanifolds, see {\bbr}. 

Our analysis focuses on the investigation of the geometry of graded
principal bundles, as well as on the reformulation of the notion of 
graded connection in terms of graded distributions and Lie
superalgebra-valued graded differential forms. We believe that our
approach makes the theory of graded connections more accessible and
close to applications in supersymmetric Yang-Mills theory, providing
thus the possibility to this elegant and economical theory to be
applied on specific problems of supersymmetric physics.   

Throughout the paper, we follow the original terminology of 
Kostant, {\kost} and we avoid the use of the term 
``supermanifold" or ``Lie supergroup". This strict distinction 
between the terms ``supermanifold"  and ``graded manifold" is 
justified if one wishes to avoid confusion with the DeWitt's or 
Roger's supermanifold theory, {\bruzzo}, {\arcad}, {\roger},
{\rogers}, {\dewitt}. In this latter approach, one constructs the
supermanifold following the pattern of ordinary differential geometry 
(using supercharts, superdifferentiability, etc.). In this way, the
supermanifold is a topological space, and in particular, a Banach
manifold, {\bruzzo}. In Kostant-Leites approach, the supermanifold is
a pair which consists of a usual differentiable manifold equipped with
a sheaf satisfying certain properties. 

Despite the fact that Kostant's graded manifolds seem to be more
abstract, we have chosen to follow Kostant's formalism for several
reasons. The most important is that in this formalism we avoid the
complications of DeWitt's supergeometry coming from
superdifferentiability of mappings. Furthermore, Kostant's graded
manifold theory is closer to ordinary manifold theory, since it
constitutes a special kind of sheaf theory over ordinary manifolds.
Finally, the algebraic concepts entering in graded manifold theory
often significantly simplify the calculus on graded manifolds. 

The article is organized as follows. In section 2, we recall 
the fundamentals of graded manifold theory, {\kost}, which 
will be necessary for the subsequent analysis, introducing at 
the same time some new concepts and tools, useful in computing
pull-backs of graded differential forms. Next section deals
with graded Lie theory. After a review of the basics of the theory,
{\kost}, {\alm}, we show how one can obtain a right action of a
graded Lie group from a left one (and vice versa); we also explain  
how a graded Lie group action gives rise to graded manifold morphisms
and derivations. A graded analog of the adjoint action of a Lie group
on itself is analyzed, and the parallelizability theorem for graded
Lie groups is proved. In section 4 we investigate the relation between
graded distributions and graded Lie group actions. The
free action in graded geometry is the central notion of this section for
which we provide two equivalent characterizations. The main result is
that each free action of a graded Lie group induces a regular graded 
distribution. We also study quotient graded structures, completing the
work of {\alm} on this subject. 

The notion and the geometry of graded principal bundle are analyzed
in sections 5 and 6. Here, we adopt a definition of this object
slightly different from that in {\almdyo}, completing in some sense, the
definition of the last reference. We prove that several
properties of ordinary principal bundles remain valid in the graded
context, making of course the appropriate modifications and
generalizations. For example, we prove that the orbits of the structure
graded Lie group are identical to the fibres of the projection and
that a graded principal bundle is globally trivial, if and only if
admits a global graded section. As an application of the tools 
developed until now, we provide an alternative definition of the 
graded principal
bundle, often very useful as we explain in section 6. After the
introduction and study of the main properties of Lie superalgebra-valued 
graded differential forms in section 7, we are ready to introduce the
notion of graded connection. Our definition is guided by the
geometrical structure of the graded principal bundle: the sheaf of
vertical derivations coincides with the graded distribution induced by
the action of the structure graded Lie group. We discuss two
equivalent definitions of the graded connection and we show how one
can construct a graded connection locally. As a direct application, we
establish the existence theorem for graded connections. 

The graded curvature is the subject of section 9. We prove the 
structure equation and the Bianchi identity for the curvature of a
graded connection. We show finally that this notion controls the
involutivity of the horizontal graded distribution determined by a
graded connection.

The previous analysis on the curvature shows also that the graded
connection $\biomega$ (as well as its curvature $F^{\biomega}$)
decompose as $\biomega=\biomega_{0}+\biomega_{1}$, where
$\biomega_{0}$ is an even graded differential 1-form with values in
the even part of the Lie superalgebra $\frak g$ of the structure
group, while $\biomega_{1}$ is an odd graded differential 1-form with
values in $\frak g_{1}$, the odd part of $\frak g$. Furthermore, the
restriction of $\biomega_{0}$ to the underlying differentiable
manifold (which is an ordinary principal bundle) gives rise to a usual
connection and exactly this observation suggests to interpret
$\biomega$, in physics terminology, as a supersymmetric gauge
potential incorporating the usual gauge potential $\biomega_{0}$ and
its supersymmetric partner $\biomega_{1}$.

\vskip0.2cm

{\bf Notational conventions.} For an algebra $A$, $A^{\ast}$ denotes
its full dual and $A^{\circ}$ its finite dual (in contrast to
Kostant's conventions where $A\tonos$ denotes the full, and $A^{\ast}$
the finite dual). If $\A$ is a sheaf of algebras over a differentiable
manifold $M$, then $\ena_{\script A}$ denotes the unit of $\A(M)$ and
$m_{\script A}$ the algebra multiplication. Throughout this article 
the term ``graded commutative" means ``$\Z_{2}$-graded commutative", 
unless otherwise stated. If $E$ is a $\Z_{2}$-graded vector
space, then $E_{0}$ and $E_{1}$ stand for its even and odd subspaces:
$E=E_{0}\oplus E_{1}$. For an element $v\an E_{i}, i=0,1$,
$|v|=i$ denotes the $\Z_{2}$-degree of $v$. Elements belonging only to
$E_{0}$ or $E_{1}$ are called homogeneous (even or odd, respectively).

\chapter{Graded manifold theory}

In this section, we review the basic notions of graded 
manifold theory, {\kost}. We also introduce some new concepts
following the pattern of ordinary differential geometry as well as
tools that will be useful in the sequel. More precisely,
the notions of vertical and projectable derivations are
introduced; push-forward (and pull-back) of derivations under
isomorphisms of graded manifolds is defined and some useful 
formulas for the computation of pull-backs of forms are 
established. We begin with the notion of the graded manifold, {\kost}.

A ringed space $(M,\A)$ is called a graded manifold of dimension
$(m,n)$ if $M$ is a differentiable manifold of dimension $m$, $\A$ is 
a sheaf of graded commutative algebras, there exist a sheaf
epimorphism $\varrho\Colon\A\rightarrow C^{\infty}$ and for each open
$U$, an open covering $\{V_{i}\}$ such that $\A(V_{i})$
is isomorphic as a graded commutative algebra to
$C^{\infty}(V_{i})\otimes\Lambda\R^{n}$, in a manner compatible with
the restriction morphisms of the sheaves involved. 
Here, $C^{\infty}$ stands for the sheaf of differentiable functions on
$M$, equipped with its trivial grading: 
$\big(C^{\infty}(U)\big)_{0}=C^{\infty}(U)$ for each open
subset $U\subset M$. We note $\dim(M,\A)=(m,n)$. 

An open $U\subset M$, for which $\A(U)\cong
C^{\infty}(U)\otimes\Lambda\R^{n}$, is called an $\A$-splitting
neighborhood. A graded coordinate system on an $\A$-splitting 
neighborhood $U$ is a collection 
$(x^{i},s^{j})=(x^{1},\ldots,x^{m},s^{1},\ldots,s^{n})$ of homogeneous
elements of $\A(U)$ with $|x^{i}|=0,|s^{j}|=1$, such that 
$(\tilde{x}^{1},\ldots,\tilde{x}^{m})$ is an ordinary coordinate
system on the open $U$, where $\tilde{x}^{i}=\varrho(x^{i})
\an C^{\infty}(U)$ and $(s^{1},\ldots,s^{n})$ are algebraically 
independent elements, that is, $\prod_{j=1}^{n}s^{j}\neq 0$.
It can also be shown, {\kost}, that if $\A^{1}(U)$ is the set of nilpotent
elements of $\A(U)$, then the sequence 
$$0\longrightarrow\A^{1}(U)\buildrel\over\longrightarrow
\A(U)\buildrel\over\longrightarrow
C^{\infty}(U)\longrightarrow 0\eqn\exactsequence$$
is exact.

Given two graded manifolds $(X,\A)$ and $(Y,\B)$, one can form their
product $(Z,\C)=(X,\A)\times(Y,\B)$ which is also a graded manifold, 
{\ruipma}, where $Z=X\times Y$ and the sheaf $\C$ is given by
$\C(Z)=\A(X)\opi\B(Y)$; in the previous tensor product, $\pi$ means
the completion of $\A(X)\otimes\B(Y)$ with respect to
Grothendieck's $\pi$-topology, {\grot}, {\ruipma}.

One can define a morphism between two graded manifolds as being a
morphism of ringed spaces compatible with the sheaf epimorphism
$\varrho$, {\ruipma}, but often it is more convenient to do this 
in a more concise way: a morphism $\sigma\Colon(M,\A)\rightarrow(N,\B)$
between two graded manifolds is just a morphism
$\sigma^{\ast}\Colon\B(N)\rightarrow\A(M)$ of graded commutative
algebras.

A very useful object in graded manifold theory is the finite dual 
$\A(M)^{\circ}$ of $\A(M)$ defined as
$$\A(M)^{\circ}=\{a\an\A(M)^{\ast}\,|\,a\,\hbox{vanishes on
an ideal of $\A(M)$ of finite codimension}\}.$$
Using general algebraic techniques (see for example {\mont}), one
readily verifies that $\A(M)^{\circ}$ is a graded cocommutative
coalgebra, the coproduct $\dacirc$ and counit $\epsilon_{\script
A}^{\circ}$ on $\A(M)^{\circ}$ being given by $\dacirc a(f\otimes
g)=a(fg)$, $\epsilon_{\script A}^{\circ}(a)=a(\ena_{\script A})$, $\forall
a\an\A(M)^{\circ}$, $f,g\an\A(M)$.
The set of group-like elements of this coalgebra contains only
elements of the form $\delta_{p}$ for $p\an M$, where
$\delta_{p}\Colon\A(M)\rightarrow\R$ is defined by
$$\delta_{p}(f)=\varrho(f)(p)=\tilde{f}(p).\eqn\deltap$$

Furthermore, if $\sigma\Colon(M,\A)\rightarrow(N,\B)$ is a
morphism of graded manifolds, then the element
$\sigma_{\ast}a\an\B(N)^{\ast}$ defined for $a\an\A(M)^{\circ}$ by
$$\sigma_{\ast}a(g)=a(\sigma^{\ast}g), 
\forall g\an\B(N)\eqn\sigmorph$$
vanishes on an ideal of finite codimension. We thus obtain a morphism
of graded coalgebras
$\sigma_{\ast}\Colon\A(M)^{\circ}\rightarrow\B(N)^{\circ}$ which
respects the group-like elements and induces a differentiable map 
$\sigma_{\ast}|_{M}\Colon M\rightarrow N$.

Another very important property of the graded coalgebra
$\A(M)^{\circ}$ is that the set of its primitive elements with respect
to the group-like element $\delta_{p}$, i.e. elements 
$v\an\A(M)^{\circ}$ for which $\dacirc(v)=v\otimes\delta_{p}+
\delta_{p}\otimes v$, is equal to the set of derivations at $p$ 
on $\A(M)$, that is, the set of elements $v\an\A(M)^{\ast}$ for which 
$$v(fg)=(vf)(\delta_{p}g)+(-1)^{|v||f|}(\delta_{p}f)(vg), \forall
f,g\an\A(M)\eqn\tangentvector$$
if $v$ and $f$ are homogeneous. This set is a subspace of
$\A(M)^{\circ}$ which we call tangent space of $(M,\A)$ at $p\an M$
and we note it by $T_{p}(M,\A)$. One easily verifies that the morphisms
of graded manifolds preserve the subspaces of primitive elements:
if $v\an T_{p}(M,\A)$, then $\sigma_{\ast}v\an T_{q}(N,\B)$
where $q=\sigma_{\ast}|_{M}(p)$. Hence, we have a well-defined notion of
the tangent (or the differential) of the morphism $\sigma$ at any
point $p\an M$ and we adopt the notation
$$T_{p}\sigma\Colon T_{p}(M,\A)\rightarrow T_{q}(N,\B),\quad
T_{p}\sigma(v)=\sigma_{\ast}v.\eqn\tangent$$  
In this context, the set $\Der\A(U)$ of derivations of $\A(U)$ plays
the r\^ole of graded vector fields on $(U,\A|_{U}), U\subset M$ open.
The difference with the ordinary differential geometry is that we
cannot evaluate directly a derivation $\xi\an\Der\A(U)$ at a point
$p\an U$ in order to obtain a tangent vector belonging to $T_{p}(M,\A)$. 
Instead, we may associate to each $\xi\an\Der\A(U)$ a tangent vector 
$\tilde{\xi}_{p}\an T_{p}(M,\A), \forall p\an U$ in the following way:
for each $f\an\A(U)$, we define $$\tilde{\xi}_{p}(f)=\delta_{p}(\xi
f).\eqn\derivation$$

If $U$ is an open subset of $M$ and $(m,n)=\dim(M,\A)$,
the set of derivations $\Der\A(U)$ is a free left $\A(U)$-module of
dimension $(m,n)$. If $U$ is an $\A$-splitting neighborhood, 
$(x^{i},s^{j})$ a graded coordinate system on $U$ and 
$\xi\an\Der\A(U)$, then there exist elements $\xi^{i},\xi^{j}\an\A(U)$
such that 
$$\xi=\sum_{i=1}^{m}\xi^{i}{\meriki\over{\meriki x^{i}}}+
\sum_{j=1}^{n}\xi^{j}{\meriki\over{\meriki s^{j}}},$$
where the derivations $\meriki/\meriki x^{i}$ and 
$\meriki/\meriki s^{j}$ are defined by
$${\meriki\over{\meriki s^{k}}}(s^{\ell})=\delta_{k}^{\ell}\ena_{U},\,
{\meriki\over{\meriki s^{k}}}(x^{i})=0,\,
{\meriki\over{\meriki x^{i}}}(s^{k})=0,\,
{\meriki\over{\meriki x^{i}}}(x^{j})=\delta_{i}^{j}\ena_{U},\eqn\defderiv$$ 
$\ena_{U}$ being the unit of $\A(U)$. Clearly, this 
decomposition is not valid in general for derivations belonging to 
$\Der\A(M)$. Therefore, we give the following definition:

\math{Definition.}{\parallel}{\sl A graded manifold $(M,\A)$ is called
parallelizable if the set of derivations $\Der\A(M)$ admits a global 
basis on $\A(M)$ consisting of $m$ even and $n$ odd derivations.}

\vskip0.3cm

The difficulty one encounters in ordinary manifold theory to
push-forward (or to pull-back) a vector field by means of a 
differentiable mapping is also present in the context of graded 
manifolds. As in the case of ordinary manifolds, this is possible 
only if we use isomorphisms of the graded manifold structure. More precisely:

\math{Definition.}{\pushfor}{\sl Let
$\sigma\Colon(M,\A)\rightarrow(N,\B)$ be an isomorphism between the
graded manifolds $(M,\A)$ and $(N,\B)$ and let $U\subset M$ be an open
subset such that $U=\sigma_{\ast}^{-1}(V), V\subset N$ open. If
$\xi\an\Der\A(U)$, then we define the push-forward
$\sigma_{\ast}\xi\an\Der\B(V)$ as
$\sigma_{\ast}\xi=(\sigma^{\ast})^{-1}\comp\xi\comp\sigma^{\ast}$.
For the pull-back of $\eta\an\Der\B(V)$, we define
$\sigma^{\ast}\eta=(\sigma^{-1})_{\ast}\eta\an\Der\A(U)$.}

\vskip0.3cm

In many situations it happens that some derivations may be related 
through a morphism (which is not necessarily an isomorphism) in a
manner similar to those of pull-back. We give the following
definition:

\math{Definition.}{\proj}{\sl Let
$\sigma\Colon(M,\A)\rightarrow(N,\B)$ be a morphism of graded
manifolds. We call two derivations $\xi\an\Der\A(M)$ et $\eta\an\Der\B(N)$
$\sigma$-related if for each $f\an\B(N)$ we have
$\sigma^{\ast}(\eta f)=\xi(\sigma^{\ast}f)$. Especially, if we fix
$\xi$ and $\sigma$ is an epimorphism, the derivation $\eta$, if it
exists, is unique;
in such a case we note $\eta=\sigma_{\ast}\xi$, we call $\xi$ a 
$\sigma_{\ast}$-projectable derivation and $\sigma_{\ast}\xi$ its 
projection by means of $\sigma$. $\xi$ will be called vertical 
derivation if it is $\sigma_{\ast}$-projectable and $\sigma_{\ast}\xi=0$.}

\vskip0.3cm

It is easy to verify that the set $\Pro(\sigma_{\ast},\A)(M)$ of
$\sigma_{\ast}$-projectable derivations is a left $\B(N)$-module and
that the set $\Ver(\sigma_{\ast},\A)(M)$
of vertical derivations is a left $\A(M)$-module. Indeed, 
for $g\an\B(N)$ and $\xi\an\Pro(\sigma_{\ast},\A)(M)$ let us
define $g\cdot\xi\an\Der\A(M)$ as $g\cdot\xi=(\sigma^{\ast}g)\xi$. Then
for $f\an\B(N)$ we find: $(g\cdot
\xi)(\sigma^{\ast}f)=(\sigma^{\ast}g)\xi(\sigma^{\ast}f)=\sigma^{\ast}g\cdot
\sigma^{\ast}[\sigma_{\ast}\xi(f)]=\sigma^{\ast}[g(\sigma_{\ast}\xi)f]$,
which proves that $g\cdot\xi\an\Pro(\sigma_{\ast},\A)(M)$ as well. We
proceed analogously for $\Ver(\sigma_{\ast},\A)(M)$. The corresponding
sheaves are denoted by
$\Pro(\sigma_{\ast},\A),\Ver(\sigma_{\ast},\A)$.

The final part of this section is devoted to a short review
of the basic properties of graded differential forms from {\kost}. 
We also develop some useful tools for computing pull-backs of 
forms using the previously mentioned concepts of push-forward 
and $\sigma$-related derivations. 

Let then $(M,\A)$ be a graded manifold of dimension $(m,n)$. For an open
$U\subset M$ we consider the tensor algebra $\yfra T(U)$ of
$\Der\A(U)$ with respect to its $\A(U)$-module structure and the ideal
$\yfra J(U)$ of $\yfra T(U)$ generated by homogeneous elements of the form 
$\xi\otimes\eta+(-1)^{|\xi||\eta|}\eta\otimes\xi$, for
$\xi,\eta\an\Der\A(U)$. Let also $\yfra T^{r}(U)\cap\yfra J(U)=
\yfra J^{r}(U)$. We call set of graded differential
$r$-forms ($r\geq 1$) on $U\subset M$ the set $\Omega^{r}(U,\A)$ of
elements belonging to ${\rm Hom}_{{\script A}(U)}(\yfra T^{r}(U),\A(U))$
which vanish on $\yfra J^{r}(U)$. For $r=0$ we define
$\Omega^{0}(U,\A)=\A(U)$ and we note the direct sum
$\displaystyle\bigoplus_{p=0}^{\infty}\Omega^{r}(U,\A)$ by 
$\Omega(U,\A)$.

If $\alpha\an\Omega^{r}(U,\A)$ and $\xi_{1},\ldots,\xi_{r}\an\Der\A(U)$,
then we denote the evaluation of $\alpha$ on the $\xi$'s by 
$(\xi_{1}\otimes\ldots\otimes\xi_{r}|\alpha)$, or simply
$(\xi_{1},\ldots,\xi_{r}|\alpha)$.

Clearly, the elements de $\Omega(U,\A)$ have a
$(\Z\oplus\Z_{2})$-bidegree and further we may define an algebra 
structure on $\Omega(U,\A)$ which thus becomes a bigraded commutative 
algebra over $\A(U)$, {\kost}. Here, we mention only the bigraded 
commutativity relation for this structure: 
if $\alpha\an\Omega^{i_{\alpha}}(U,\A)_{j_{\alpha}},
\beta\an\Omega^{i_{\beta}}(U,\A)_{j_{\beta}}$, then 
$\alpha\beta\an\Omega^{i_{\alpha}+i_{\beta}}(U,\A)_{j_{\alpha}+j_{\beta}}$
and $\alpha\beta=(-1)^{i_{\alpha}i_{\beta}+j_{\alpha}j_{\beta}}
\beta\alpha$.

Next consider the linear map
$d\Colon\Omega^{0}(U,\A)\rightarrow\Omega^{1}(U,\A)$ defined by
$$(\xi|dg)=\xi(g), \,\xi\an\Der\A(U), g\an\A(U).\eqn\orismosd$$
In a graded coordinate system $(x^{i},s^{j}), i=1,\ldots
m,j=1,\ldots n$ on $U$, we take
$$dg=\sum_{i=1}^{m}dx^{i}{\meriki g\over\meriki
x^{i}}+\sum_{j=1}^{n}ds^{j}{\meriki g\over\meriki s^{j}}.\eqn\dfison$$
One can extend this linear map to a derivation
$d\Colon\Omega(U,\A)\rightarrow\Omega(U,\A)$ 
of bidegree (1,0) such that $d^{2}=0$ and $d|_{\Omega^{0}(U,\script A)}$
gives equation $\dfison$. We will call $d$ exterior differential on graded 
differential forms. 

Interior products and Lie derivatives with
respect to elements of $\Der\A(U)$ also make sense in the graded setting. 
Indeed, if $\alpha\an\Omega^{r+1}(U,\A)$, then for
$\xi,\xi_{1},\ldots,\xi_{r}$ homogeneous elements of $\Der\A(U)$
we define
$$(\xi_{1},\ldots,\xi_{r}|{\bit
i}(\xi)\alpha)=(-1)^{|\xi|\sum_{i=1}^{r}|\xi_{i}|}
(\xi,\xi_{1},\ldots,\xi_{r}|\alpha)\eqn\interiorproduct$$
and we thus obtain a linear map 
${\bit i}(\xi)\Colon\Omega(U,\A)\rightarrow\Omega(U,\A)$
of bidegree $(-1,|\xi|)$. This is the interior product with respect
to $\xi$. 

Lie derivatives are defined  as usual by means of Cartan's algebraic
formula: $\Lie_{\xi}=d\comp{\bit i}(\xi)+{\bit i}(\xi)\comp d$, 
thus $\Lie_{\xi}$ has bidegree $(0,|\xi|)$. Furthermore, it can be
proved that the morphism of graded commutative algebras
$\sigma^{\ast}\Colon\B(W)\rightarrow\A(U),
U=\sigma^{-1}_{\ast}(W)\subset M$ coming from a morphism 
$\sigma\Colon(M,\A)\rightarrow(N,\B)$ of graded manifolds,
can be extended to a unique morphism of bigraded commutative algebras
$\sigma^{\ast}\Colon\Omega(W,\B)\rightarrow\Omega(U,\A)$, which
commutes with the exterior differential.

As for the case of derivations, see relation $\derivation$, one can define
for each graded differential form $\alpha\an\Omega^{r}(U,\A)$ a
multilinear form $\tilde{\alpha}_{p}$ (with real values) on the 
tangent space $T_{p}(M,\A)$ for each $p\an U$. It suffices to set 
$$\big((\tilde{\xi}_{1})_{p},\ldots,
(\tilde{\xi}_{r})_{p}|\tilde{\alpha}_{p}\big)=\delta_{p}
(\xi_{1},\ldots,\xi_{r}|\alpha).\eqn\multiliforms$$
The set of forms on $U$ obtained in this way is denoted by
$\Omega^{r}_{\script A}(U)$.

We establish now a general method for the calculation of pull-backs of
graded differential forms.

\math{Proposition.}{\pullback}{\sl Let $\sigma\Colon(M,\A)\rightarrow
(N,\B)$ be a morphism of graded manifolds, $W\subset N$ and $U=
\sigma^{-1}_{\ast}(W)\subset M$. Then, if $\alpha\an\Omega^{r}(W,\B)$ and
$\xi_{i}\an\Der\A(U)$ and $\eta_{i}\an\Der\B(W)$ are $\sigma$-related
for $i=1,\ldots,r$, we have:
$$(\xi_{1},\ldots,\xi_{r}|\sigma^{\ast}\alpha)=
\sigma^{\ast}(\eta_{1},\ldots,\eta_{r}|\alpha).\eqn\pullba$$
In particular, when $\sigma$ is an isomorphism, the previous relation
holds for each $\xi_{i}\an\Der\A(U)$, setting $\eta_{i}=\sigma_{\ast}
\xi_{i}$.}

\indent \undertext{\it Proof}. Using the fact that $\sigma^{\ast}$ is an
isomorphism of bigraded commutative algebras and that $\Omega(W,\B)$
is the exterior algebra of $\Omega^{1}(W,\B)$, it is sufficient to
prove this formula for $\alpha=df, f\an\B(W)$. If $\xi\an\Der\A(U)$
and $\eta\an\Der\B(W)$ are $\sigma$-related, then 
$(\xi|\sigma^{\ast}\alpha)=(\xi|d(\sigma^{\ast}f))=
\xi(\sigma^{\ast}f)=\sigma^{\ast}(\eta
f)=\sigma^{\ast}(\eta|df)=\sigma^{\ast}(\eta|\alpha)$. In particular,
when $\sigma$ is an isomorphism, each $\xi$ is $\sigma$-related to
$\eta=\sigma_{\ast}\xi$.\qed

\vskip0.3cm

One easily proves that pull-backs on forms belonging to
$\Omega^{r}_{\script B}(W)$ can be expressed using the 
familiar formulas of ordinary differential geometry but taking as
tangent of the morphism $\sigma$, the linear mapping defined through
relations $\sigmorph$, $\tangentvector$.

\chapter{Elements of graded Lie theory}

As we will see later in detail, graded Lie theory plays a central 
r\^ole in the geometry of graded principal bundles just as ordinary 
theory of Lie groups does in differentiable principal bundles. This 
section deals with the notion and elementary properties of graded 
Lie groups; more information can be found in {\kost}. Furthermore, 
we prove some facts about the theory of actions of graded Lie groups, 
which, to the best of our knowledge, have not ever appeared in the 
literature. We first give the definition of a graded Lie group,
{\alm}, {\kost}. 

\math{Definition.}{\gradedliegroup}{\sl A graded Lie group $(G,\A)$
is a graded manifold such that $G$ is an ordinary Lie group, the
algebra $\A(G)$ is equipped with the structure of a graded 
Hopf algebra with antipode and furthermore, the algebra epimorphism 
$\varrho\Colon\A(G)\rightarrow C^{\infty}(G)$ is a morphism of 
graded Hopf algebras.}

\math{Remark.}{\completionpi}{\it In the graded Hopf algebra structure of
the previous definition, all tensor products are completions of the
usual ones with respect to Grothendieck's $\pi$-topology, {\rm{\alm}},
{\rm{\ruipma}}.} 

\vskip0.2cm

We denote by $\da$, $\epsilon_{\script A}$, $s_{\script A}$ the
coproduct, counit and antipode of $\A(G)$, respectively.
It is possible to prove that the finite dual $\A(G)^{\circ}$
inherits also a graded Hopf algebra structure from $\A(G)$.
The algebra multiplication on $\A(G)^{\circ}$ is given by
the convolution product:
$$(a\odot b)=(a\otimes b)\comp\da, \,\forall a,b\an\A(G)^{\circ},
\eqn\convolution$$
and the unit of $\A(G)^{\circ}$ with respect to $\odot$ is 
the counit of $\A(G)$. One proves easily that the set of primitive
elements of $\A(G)^{\circ}$ with respect to $\delta_{e}$ ($e$ is the
identity of $G$), that is, the tangent space $T_{e}(G,\A)$, is a graded Lie
algebra, the bracket being given by $[u,v]=u\odot
v-(-1)^{|u||v|}v\odot u$, for homogeneous elements $u,v$. We call
$T_{e}(G,\A)$ Lie superalgebra of $(G.\A)$ and denote it by $\frak
g$. Clearly, $\frak g=\frak g_{0}\oplus\frak g_{1}$, where $\frak
g_{0}=T_{e}G$.

Using also the fact that $\varrho\Colon\A(G)\rightarrow C^{\infty}(G)$
is a morphism of graded Hopf algebras one readily verifies that the
convolution product $\odot$ is compatible with the group structure of
$G$ in the sense that $\delta_{g}\odot\delta_{h}=\delta_{gh}$. 

A very important property of the finite dual $\A(G)^{\circ}$ is given
by the following, {\kost}.

\math{Theorem.}{\liehopfalgebra}{\sl For each graded Lie group
$(G,\A)$, the finite dual $\A(G)^{\circ}$ has the structure of a 
Lie-Hopf algebra. In fact, $\A(G)^{\circ}=\R(G)\crossed E(\frak g)$,
where $\R(G)$ is the group algebra of $G$, $E(\frak g)$ is the
universal enveloping algebra of $\frak g$ and $\crossed$ is the smash
product of $\R(G)$ and $E(\frak g)$ with respect to the adjoint
representation of $G$ on the superalgebra $\frak g$.}

\vskip0.3cm

For the adjoint representation of $G$ on $\frak g$, see also below in
this section. 

There exist, in this setting, analogs of left and right translations
on a Lie group.

\math{Definition.}{\leftrighttrans}{\sl Let $(G,\A)$ be a graded Lie
group and $a\an\A(G)^{\circ}$. Set $r_{a}=(id\otimes a)\comp\da$
and $\ell_{a}=(a\otimes id)\comp\da$. We call the endomorphisms $r_{a}$
and $\ell_{a}$ of $\A(G)$ right and left translations respectively on
$(G,\A)$.}

\vskip0.3cm

In order to justify this terminology, we first discuss some
important properties of these endomorphisms.

\math{Proposition.}{\leftrightprop}{\sl 
\item{1.} $r_{\epsilon_{\script A}}=\ell_{\epsilon_{\script A}}=id$
\item{2.} $r_{a\odot b}=r_{a}\comp r_{b}$, 
$\ell_{a\odot b}=(-1)^{|a||b|}\ell_{b}\comp\ell_{a}$
\item{3.} $r_{b}\comp\ell_{a}=(-1)^{|a||b|}\ell_{a}\comp r_{b}, \forall
a,b\an\A(G)^{\circ}$
\item{4.} If $a\an\A(G)^{\circ}$ is group-like, then $r_{a}$ and
$\ell_{a}$ are graded algebra isomorphisms.

}

\vskip0.3cm

We postpone the proof until the study of actions of graded Lie groups 
(see below in this section). Then, it will be clear that the previous
proposition is immediate using the general techniques of actions (note
that Proposition {\leftrightprop} has already appeared in {\alm} without
proof).

Part (4) of this proposition tells us that if $a$ is group-like
$a=\delta_{g}$, $g\an G$, then there exist morphisms of graded 
manifolds $R_{g}\Colon(G,\A)\rightarrow(G,\A)$, 
$L_{g}\Colon(G,\A)\rightarrow(G,\A)$ such that 
$R^{\ast}_{g}=r_{a}$ and $L^{\ast}_{g}=\ell_{a}$. It is interesting
to calculate the coalgebra morphisms 
$R_{g\ast}\Colon\A(G)^{\circ}\rightarrow\A(G)^{\circ}$ and 
$L_{g\ast}\Colon\A(G)^{\circ}\rightarrow\A(G)^{\circ}$.
Consider for example $R_{g\ast}$. If $b\an\A(G)^{\circ}$ and
$f\an\A(G)$, we find:
$$R_{g\ast}b(f)=b(r_{a}f)=b\big(\sum_{i}
(-1)^{|a||I^{i}f|}(I^{i}f) a(J^{i}f)\big)
=\sum_{i}(b\otimes a)(I^{i}f\otimes J^{i}f),$$
where we have set $\da f=\sum_{i}I^{i}f\otimes J^{i}f$. Hence,
$R_{g\ast}b=b\odot a$ and similarly $L_{g\ast}b=a\odot b$. 
This means that $r_{a}$ and $\ell_{a}$ correspond to right and left
translations, as one can see at the coalgebra level. 

Next, we introduce the graded analog of actions (see also {\alm}). 
Let then $(G,\A)$ be a graded Lie group and $(Y,\B)$ a graded
manifold. We give the following definition.

\math{Definition.}{\action}{\sl We say that the graded Lie group
$(G,\A)$ acts on the graded manifold $(Y,\B)$ to the right if there exists a 
morphism $\Phi\Colon(Y,\B)\times(G,\A)\rightarrow(Y,\B)$ of graded
manifolds such that the corresponding morphism of graded commutative
algebras $\Phi^{\ast}\Colon\B(Y)\rightarrow\B(Y)\opi\A(G)$ defines
a structure of right $\A(G)$-comodule on $\B(Y)$. Using the notion 
of left comodule, we may define the left action of $(G,\A)$ on $(Y,\B)$.}

\vskip0.3cm

More explicitly, if $\Phi$ is a right and $\Psi$ is a left action, 
then the morphisms
$\Phi^{\ast}$, $\Psi^{\ast}$  satisfy the following properties:
$$(id\otimes\da)\comp\Phi^{\ast}=(\Phi^{\ast}\otimes id)\comp\Phi^{\ast},
\,\,(id\otimes\epsilon_{\script
A})\comp\Phi^{\ast}=id,\eqn\rightcomodulestructure$$ 
$$(\da\otimes
id)\comp\Psi^{\ast}=(id\otimes\Psi^{\ast})\comp\Psi^{\ast},
\,\,(\epsilon_{\script A}\otimes
id)\comp\Psi^{\ast}=id.\eqn\leftcomodulestructure$$

Let now $\Psi^{r}\Colon(Y,\B)\times(G,\A)\rightarrow(Y,\B)$ be the
morphism of graded manifolds defined by 
$$\Psi^{r\ast}=(id\otimes s_{\script A})\comp
T\comp\Psi^{\ast},\eqn\psirightactiondef$$
where $T$ is the twist morphism, $T(a\otimes b)=(-1)^{|a||b|}b\otimes
a$. 

\math{Lemma.}{\psirightactionlemma}{\sl The morphism $\Psi^{r}$
defined by $\psirightactiondef$ is a right action of $(G,\A)$ on
$(Y,\B)$.}

\undertext{\it Proof}. It suffices to prove that relations
$\rightcomodulestructure$ are valid for the morphism $\Psi^{r\ast}$. 
We take: 
$(id\otimes\da)\comp\Psi^{r\ast}=
\big(id\otimes(s_{\script
A}\otimes s_{\script A})\comp T\comp\da\big)\comp T\comp\Psi^{\ast}=
\big(id\otimes (s_{\script A}\otimes s_{\script A})\comp T\big)\comp
T\comp(\da\otimes id)\comp\Psi^{\ast}=
\big(id\otimes (s_{\script A}\otimes s_{\script A})\comp T\big)\comp
T\comp\break(id\otimes\Psi^{\ast})\comp\Psi^{\ast}=
(id\otimes s_{\script A}\otimes s_{\script
A})\comp(T\comp\Psi^{\ast}\otimes id)\comp
T\comp\Psi^{\ast}=(\Psi^{r\ast}\otimes id)\comp\Psi^{r\ast}$.
On the other hand,
$$(id\otimes\epsilon_{\script A})\comp\Psi^{r\ast}=
(id\otimes\epsilon_{\script A}\comp s_{\script A})\comp
T\comp\Psi^{\ast}=(id\otimes\epsilon_{\script A})\comp T\comp\Psi^{\ast}=id,$$
which completes the proof.\qed

\vskip0.3cm

Similarly, for a right action $\Phi\Colon(Y,\B)\times(G,\A)
\rightarrow(Y,\B)$ one can define in a
canonical way, a left action
$\Phi^{\ell}\Colon(G,\A)\times(Y,\B)\rightarrow(Y,\B)$ as
$\Phi^{\ell\ast}=(s_{\script A}\otimes id)\comp T\comp\Phi^{\ast}$. 
Then, the restriction $\Phi_{\ast}|_{Y\times G}\Colon Y\times
G\rightarrow Y$ defines a right action of $G$ on the manifold $Y$ and
furthermore, for the canonically associated left action $\Phi^{\ell}$,
the restriction $\Phi^{\ell}_{\ast}|_{G\times Y}\Colon G\times
Y\rightarrow Y$ is a left action given by
$\Phi^{\ell}_{\ast}|_{G\times Y}(g,y)=\Phi_{\ast}|_{Y\times
G}(y,g^{-1})$, as one expects. We have analogous facts for the left
action $\Psi$. Observe here that the possibility to define
$\Phi^{\ell\ast}$ and $\Psi^{r\ast}$ as morphisms of graded
commutative algebras depends crucially on the fact that the antipode
$s_{\script A}\Colon\A(G)\rightarrow\A(G)$ is a morphism of graded
commutative algebras.

\math{Remark.}{\paratirisidyo}{\it The right action $\Phi^{\ell r}$
canonically associated to $\Phi^{\ell}$ equals to $\Phi$:
$\Phi^{\ell r\ast}=(id\otimes s_{\script A})\comp T\comp(s_{\script
A}\otimes id)\comp T\comp\Phi^{\ast}=(id\otimes s_{\script A})\comp
T^{2}\comp(id\otimes s_{\script A})\comp\Phi^{\ast}=\Phi^{\ast}$,
since $T^{2}=id, s_{\script A}^{2}=id$.}

\vskip0.2cm

For a right action $\Phi$, one may introduce for each
$a\an\A(G)^{\circ}$ and $b\an\B(Y)^{\circ}$, two linear
maps $(\Phi^{\ast})_{a}\Colon\B(Y)\rightarrow\B(Y)$ and 
$(\Phi^{\ast})_{b}\Colon\B(Y)\rightarrow\A(G)$ as follows:
$$(\Phi^{\ast})_{a}=(id\otimes a)\comp\Phi^{\ast}\quad\hbox{and}\quad 
(\Phi^{\ast})_{b}=(b\otimes id)\comp\Phi^{\ast}.\eqn\linearmapsena$$
Similarly, for a left action $\Psi$ one defines
$$(\Psi^{\ast})_{a}=(a\otimes id)\comp\Psi^{\ast}\quad\hbox{and}\quad 
(\Psi^{\ast})_{b}=(id\otimes b)\comp\Psi^{\ast}.\eqn\linearmapsdyo$$
The following theorem clarifies the r\^ole of these maps.

\math{Theorem.}{\thewrima}{\sl 
\item{1.} $(\Phi^{\ast})_{\epsilon_{\script
A}}=(\Psi^{\ast})_{\epsilon_{\script A}}=id$
\item{2.} $(\Phi^{\ast})_{a_{1}\odot
a_{2}}=(\Phi^{\ast})_{a_{1}}\comp(\Phi^{\ast})_{a_{2}}$, 
$(\Psi^{\ast})_{a_{1}\odot
a_{2}}=(-1)^{|a_{1}||a_{2}|}
(\Psi^{\ast})_{a_{2}}\comp(\Psi^{\ast})_{a_{1}}$
\item{3.}$(\Phi^{\ast})_{b}\comp(\Phi^{\ast})_{a}=
(-1)^{|a||b|}r_{a}\comp(\Phi^{\ast})_{b}$, 
$(\Psi^{\ast})_{b}\comp(\Psi^{\ast})_{a}=
(-1)^{|a||b|}\ell_{a}\comp(\Psi^{\ast})_{b}$
\item{4.} If $a,b$ are group-like elements, then $(\Phi^{\ast})_{a}$
is an isomorphism and $(\Phi^{\ast})_{b}$ is a morphism of graded
commutative algebras. In particular, if $a=\delta_{g},b=\delta_{y}$,
then we write the corresponding morphisms of graded manifolds as
$\Phi_{g}\Colon\break(Y,\B)\rightarrow(Y,\B)$ and
$\Phi_{y}\Colon(G,\A)\rightarrow(Y,\B)$, so
$\Phi_{g}^{\ast}=(\Phi^{\ast})_{\delta_{g}}$ and 
$\Phi_{y}^{\ast}=(\Phi^{\ast})_{\delta_{y}}$. Similarly for
$(\Psi^{\ast})_{a},(\Psi^{\ast})_{b}$. 
\item{5.} If $a$ is a primitive element with respect to $\delta_{e}$,
then $(\Phi^{\ast})_{a},(\Psi^{\ast})_{a}\an\Der\B(Y)$. We call these 
derivations the induced (by the action and the element $a$) 
derivations on $\B(Y)$.

}

\undertext{\it Proof}.
 
(1) Evident, by the defining properties of the left and right action.

(2) We prove this property only for the right action $\Phi$; one
proceeds in a similar way for the left action $\Psi$. We have:
$$\eqalign{(\Phi^{\ast})_{a_{1}\odot a_{2}}&=(id\otimes(a_{1}\odot
a_{2}))\comp\Phi^{\ast}=(id\otimes a_{1}\otimes
a_{2})\comp(id\otimes\da)\comp\Phi^{\ast}\cr
\hfill&=(id\otimes a_{1}\otimes
a_{2})\comp(\Phi^{\ast}\otimes id)\comp\Phi^{\ast}
=(id\otimes a_{1})\comp\Phi^{\ast}\comp(id\otimes
a_{2})\comp\Phi^{\ast}\cr
\hfill&=(\Phi^{\ast})_{a_{1}}\comp(\Phi^{\ast})_{a_{2}}.\cr}$$

(3) Again, we give the proof only for the right action.
$$\eqalign{(\Phi^{\ast})_{b}\comp(\Phi^{\ast})_{a}&=
(b\otimes id)\comp\Phi^{\ast}\comp(id\otimes a)\comp\Phi^{\ast}\cr
\hfill&=(-1)^{|a||b|}(id\otimes a)(b\otimes
id)\comp(id\otimes\da)\comp\Phi^{\ast}\cr
\hfill&=(-1)^{|a||b|}(id\otimes a)\comp\da\comp(b\otimes
id)\comp\Phi^{\ast}=(-1)^{|a||b|}r_{a}\comp(\Phi^{\ast})_{b}.\cr}$$

(4) The fact that these maps are morphisms is evident because
they are compositions of morphisms when $a,b$ are group-like.
Furthermore, $(\Phi^{\ast})_{a}$ and $(\Psi^{\ast})_{a}$
are isomorphisms because their inverses exist, as one can check 
from parts (1) and (2).

(5) A derivation $\xi$ on the graded commutative algebra $\B(Y)$ has
the property $\xi(fg)=\xi(f)g+(-1)^{|\xi||f|}f\xi(g)$ for each
homogeneous element $f\an\B(Y)$. This can be restated as follows:
$\xi\comp m_{\script B}=m_{\script B}\comp(\xi\otimes
id+id\otimes\xi)$. We call the endomorphism $\xi$ primitive, so
the set of derivations coincides with the set of primitive
elements. Using this terminology, we must prove that
$(\Phi^{\ast})_{a}$ and $(\Psi^{\ast})_{a}$ are primitive when 
$a$ is primitive with respect to $\delta_{e}$. Consider for example
$\xi=(\Phi^{\ast})_{a}$. We have:
$$\eqalign{\xi\comp m_{\script B}&
=(m_{\script B}\otimes a\comp m_{\script A})\comp(id\otimes
T\otimes id)\comp(\Phi^{\ast}\otimes\Phi^{\ast})\cr
\hfill&=(m_{\script B}\otimes(a\otimes\delta_{e}+\delta_{e}\otimes
a))\comp(id\otimes T\otimes id)\comp(\Phi^{\ast}\otimes\Phi^{\ast})\cr
\hfill&=m_{\script B}\comp[(id\otimes
a)\comp\Phi^{\ast}\otimes id]+m_{\script B}\comp
[id\otimes(id\otimes a)\comp\Phi^{\ast}]\cr
\hfill&=m_{\script B}\comp(\xi\otimes id+id\otimes\xi).\cr}$$
One proceeds similarly for $(\Psi^{\ast})_{a}$. We note finally that
if $a$ is homogeneous, then
$|(\Phi^{\ast})_{a}|=|(\Psi^{\ast})_{a}|=|a|$.\qed

\vskip0.2cm

\math{Corollary.}{\porisma}{\sl Proposition {\leftrightprop}.}

\undertext{\it Proof}. Since the coproduct $\da$ on the Hopf algebra
$\A(G)$ has the properties $(id\otimes\da)\comp\da=(\da\otimes
id)\comp\da$ and $(id\otimes\epsilon_{\script
A})\comp\da=(\epsilon_{\script A}\otimes id)\comp\da=id$, it defines
left and right actions $L$ and $R$ respectively of $(G,\A)$ on
itself. Choosing thus $(Y,\B)=(G,\A)$ in the previous theorem, we may
write $(L^{\ast})_{a}=\ell_{a}$, $(R^{\ast})_{a}=r_{a}$; Proposition 
{\leftrightprop} is then immediate.\qed

\vskip0.3cm

Next, let $(G,\A)$ be a graded Lie group and
$\ell\Colon\A(G)\rightarrow\A(G)\opi\A(G)$
the linear map defined by $$\ell=[m_{\script A}\comp(id\otimes
s_{\script A})\otimes id]\comp(id\otimes
T)\comp(id\otimes\da)\comp\da.\eqn\linearmap$$

\math{Proposition.}{\adjointaction}{\sl The linear map $\ell$ is a
morphism of graded commutative algebras defining thus a morphism of 
graded manifolds which we denote by
$AD\Colon(G,\A)\times(G,\A)\rightarrow(G,\A)$. Furthermore, 
$AD$ is a left action of $(G,\A)$ on itself.}

\undertext{\it Proof}. $\ell$ is a morphism of graded algebras as
composition of morphisms; so we can write 
$\ell=AD^{\ast}$ for a morphism of graded manifolds
$AD\Colon(G,\A)\times(G,\A)\rightarrow(G,\A)$. 
We now check relations $\leftcomodulestructure$ for $AD^{\ast}$. 
For the first one, the following identity is the key of
the proof: $(id\otimes id\otimes id\otimes \da)\comp(id\otimes
id\otimes\da)\comp(\da\otimes id)\comp\break\da=(id\otimes
id\otimes\da\otimes id)\comp(id\otimes\da\otimes id)
\comp(id\otimes\da)\comp\da$. Indeed, applying the two members of the
previous identity on the same $f\an\A(G)$, we find after a long and
cumbersome calculation the first of $\leftcomodulestructure$.
For the second of $\leftcomodulestructure$, we proceed as follows:
$$\eqalign{(\epsilon_{\script A}\otimes id)\comp\ell&=
(\epsilon_{\script A}\comp m_{\script A}\otimes
id)\comp(id\otimes s_{\script A}\otimes id)\comp(id\otimes
T)\comp(id\otimes\da)\comp\da\cr
\hfill&=(\epsilon_{\script A}\otimes\epsilon_{\script A}\comp
s_{\script A}\otimes id)\comp(id\otimes T)
\comp(id\otimes\da)\comp\da\cr
\hfill&=(\epsilon_{\script A}\otimes id\otimes\epsilon_{\script
A})\comp(id\otimes\da)\comp\da=(\epsilon_{\script A}\otimes
id)\comp\da=id,\cr}$$
which completes the proof.\qed

\vskip0.3cm

We call the action $AD$ adjoint action of $(G,\A)$ on itself. As in
ordinary Lie theory, the adjoint action respects the primitive
elements with respect to $\delta_{e}$ in the sense of the following
proposition. 

\math{Proposition.}{\adjointactionprop}{\sl Let $AD_{\ast
a}\Colon\A(G)^{\circ}\rightarrow\A(G)^{\circ}, a\an\A(G)^{\circ}$ be
defined as $AD_{\ast a}(b)=AD_{\ast}(a\otimes b)$. Then, for an
element $a\an\A(G)^{\circ}$ group-like or primitive with respect to 
$\delta_{e}$, $AD_{\ast a}$ is a linear map on the Lie superalgebra 
$\frak g$.} 

\undertext{\it Proof}. Consider first the case where $a$ is a
group-like element, $a=\delta_{g},g\an G$. If $v\an\frak g$, we have:
$AD_{\ast a}(v)=AD_{\ast}(a\otimes v)=a\odot v\odot a^{-1}$, because for
group-like elements the antipode $s_{\script A}^{\circ}$ of
$\A(G)^{\circ}$ is given by $s_{\script
A}^{\circ}a=a^{-1}=\delta_{g^{-1}}$. It is then immediate to verify
that if $\dacirc$ is the coproduct of $\A(G)^{\circ}$, we have:
$\dacirc(AD_{\ast a}(v))=AD_{\ast
a}(v)\otimes\delta_{e}+\delta_{e}\otimes AD_{\ast a}(v)$ which means
that $AD_{\ast a}(v)$ belongs also to $\frak g$. Proceeding in the
same way for the case where $a$ is primitive with respect to
$\delta_{e}$, we find $AD_{\ast a}(v)=a\odot v-(-1)^{|a||v|}v\odot
a=[a,v]\an\frak g$. Thus, for an element $a\an\A(G)^{\circ}$
group-like or primitive with respect to $\delta_{e}$, we take 
$AD_{\ast a}\an\Endo\frak g$.\qed

\vskip0.3cm

The previous proof makes clear that if $a=\delta_{g}$, $g\an G$,
then $AD_{\ast a}$ is an isomorphism of the Lie 
superalgebra $\frak g$. Indeed, in this case
$AD_{\ast a}=R_{g\ast}^{-1}\comp L_{g\ast}$. For
the case where $a$ is primitive with respect to $\delta_{e}$,
$AD_{\ast a}$ coincides with the adjoint representation of $\frak
g$ on itself, $AD_{\ast a}=ad(a)$, where $ad(a)(b)=[a,b],\forall
b\an\frak g$.

\math{Remark.}{\epeksigisidyo}{\it One can define a linear map 
$\Psi_{\ast a}\Colon\B(Y)^{\circ}\rightarrow\B(Y)^{\circ}$ 
for each left action $\Psi\Colon(G,\A)\times(Y,\B)\rightarrow(Y,\B)$
and $a\an\A(G)^{\circ}$ by $\Psi_{\ast a}(b)=\Psi_{\ast}(a\otimes b)$,
$\forall b\an\B(Y)^{\circ}$. It is then easily verified that for a
group-like, $\Psi_{\ast a}$ is an isomorphism of graded coalgebras;
furthermore, $\Psi_{\ast a}=\Psi_{g\ast}$ if $a=\delta_{g}$ (see
Theorem {\thewrima}). We have analogous facts for a right action.}

\vskip0.2cm

Graded Lie groups provide an important and wide class of
parallelizable graded manifolds as the following theorem asserts.

\math{Theorem.}{\paralleltheorem}{\sl Each graded Lie group $(G,\A)$
is a parallelizable graded manifold.}

\undertext{\it Proof}. Let
$\Phi\Colon(G,\A)\times(G,\A)\rightarrow(G,\A)$ be the right action
such that $\Phi^{\ast}=\da$ (see the proof of Corollary {\porisma}).
Then, $(\Phi^{\ast})_{a}=(id\otimes a)\comp\Phi^{\ast}=r_{a}$. If 
$a\an\frak g$, then we know by Theorem {\thewrima} that
$(\Phi^{\ast})_{a}$ is a derivation on $\A(G)$. For $g\an G$,
the tangent vector $\widetilde{(\Phi^{\ast})_{a}}(g)\an T_{g}(G,\A)$ is
calculated by means of $\derivation$: 
$\widetilde{(\Phi^{\ast})_{a}}(g)=\delta_{g}\comp(\Phi^{\ast})_{a}
=\delta_{g}\odot a=L_{g\ast}(a)$; but $L_{g\ast}$ is an isomorphism 
by Proposition {\leftrightprop}. 
This means that if $\{a^{i},b^{j}\}$ 
is a basis of the Lie superalgebra $\frak g$, $a^{i}\an\frak
g_{0},b^{j}\an\frak g_{1}$, then
$\left\{\widetilde{(\Phi^{\ast})_{a^{i}}}(g),
\widetilde{(\Phi^{\ast})_{b^{j}}}(g)\right\}$
is a basis of $T_{g}(G,\A)$ for each $g\an G$. By Proposition 2.12.1 of
{\kost}, we conclude that
$\{(\Phi^{\ast})_{a^{i}},(\Phi^{\ast})_{b^{j}}\}$ is a global basis of 
$\Der\A(G)$ for its left $\A(G)$-module structure.\qed

\chapter{Actions, graded distributions and quotient structures}

Two important notions in the study of graded Lie group actions on
graded manifolds are those of graded distributions and quotient graded
manifolds. In this section, we investigate the relation between graded
distributions and free actions properly defined in the graded setting.
In addition, we find a necessary and sufficient condition in order
that the quotient defined by the action of a graded Lie group be a
graded manifold. 

We first introduce the notion of the graded distribution (see also
{\riccigia}).

\math{Definition.}{\gdistribution}{\sl Let $(M,\A)$ be a graded
manifold of dimension $(m,n)$. We call graded distribution of
dimension $(p,q)$ on $(M,\A)$ a sheaf $U\rightarrow\E(U)$ of free
$\A(U)$-modules such that each $\E(U)$ is a graded submodule 
of $\Der\A(U)$ of dimension $(p,q)$. The distribution will be called
involutive if for each $\xi,\eta\an\E(U)$, we have $[\xi,\eta]\an\E(U),
\forall U\subset M$ open.}

\vskip0.3cm

Thus, for each open $U\subset M$, there exist elements
$\xi_{i}\an\big(\Der\A(U)\big)_{0},i=1,\ldots,p$,
$\eta_{j}\an\big(\Der\A(U)\big)_{1},j=1,\ldots,q$ such that
$$\E(U)=\A(U)\cdot\xi_{1}\oplus\cdots\oplus\A(U)\cdot\xi_{p}
\oplus\A(U)\cdot\eta_{1}\oplus\cdots\oplus\A(U)\cdot\eta_{q}.$$
Given the graded distribution $\E$ on $(M,\A)$, one obtains, for each
$x\an M$, a graded subspace $E_{x}$ of $T_{x}(M,\A)$, calculating
the tangent vectors $\tilde{\xi}_{x}\an T_{x}(M,\A)$ via relation 
$\derivation$, for each $\xi\an\E(U)$, $x\an U$. Clearly,
$E_{x}=(E_{x})_{0}\oplus(E_{x})_{1}$ with
$\dim(E_{x})_{0}=\epsilon_{0}(x)\leq p$ and 
$\dim(E_{x})_{1}=\epsilon_{1}(x)\leq q$. Therefore, we make the
following distinction:

\math{Definition.}{\regdistribution}{\sl A graded distribution $\E$ of
dimension $(p,q)$ on the graded manifold $(M,\A)$ is called 0-regular
if $\epsilon_{0}(x)=p$, and 1-regular if $\epsilon_{1}(x)=q$, for each
$x\an M$. We say that $\E$ is regular if $\epsilon_{0}(x)=p$ and
$\epsilon_{1}(x)=q$, for each $x\an M$.}

\vskip0.3cm

For the subsequent analysis, a graded generalization of free actions
will be necessary.

\math{Definition.}{\freeaction}{\sl We call the right action
$\Phi\Colon(Y,\B)\times(G,\A)\rightarrow(Y,\B)$ free, if for each
$y\an Y$ the morphism $\Phi_{y}\Colon(G,\A)\rightarrow(Y,\B)$ is such
that $\Phi_{y\ast}\Colon\A(G)^{\circ}\rightarrow\B(Y)^{\circ}$ is
injective.}

\vskip0.3cm

Similarly, one defines the left free action. It is clear that if the
graded Lie group $(G,\A)$ acts freely on $(Y,\B)$, then we obtain a
free action of $G$ on $Y$, but if only the restriction
$\Phi_{\ast}|_{Y\times G}$ is a free action then, in general, the
action $\Phi$ is not free.

An equivalent characterization of the free action in graded Lie theory
is provided by the following proposition for the case of a right
action. 

\math{Proposition.}{\freeactionequiv}{\sl The action
$\Phi\Colon(Y,\B)\times(G,\A)\rightarrow(Y,\B)$ is free if and only if
the morphism of graded manifolds $\tilde{\Phi}=
(\Phi\times\pi_{1})\comp\Delta\Colon(Y,\B)\times(G,\A)
\rightarrow(Y,\B)\times(Y,\B)$ is such that 
$\tilde{\Phi}_{\ast}$ is injective. 
Here, $\Delta$ denotes the diagonal morphism on $(Y,\B)\times(G,\A)$ 
and $\pi_{1}$ is the projection on the first factor.}

\undertext{\it Proof}. Consider elements 
$a=\delta_{g}\an\A(G)^{\circ}, b=\delta_{y}\an\B(Y)^{\circ}$
group-like and $u\an T_{y}(Y,\B), w\an T_{g}(G,\A)$ primitive. 
Then, a simple calculation gives 
$$\tilde{\Phi}_{\ast}(b\otimes a)=\Phi_{y\ast}(a)\otimes
b\eqn\relena$$  
$$\tilde{\Phi}_{\ast}(u\otimes a+b\otimes w)=
[\Phi_{g\ast}(u)+\Phi_{y\ast}(w)]\otimes b+
\Phi_{g\ast}(b)\otimes u.\eqn\reldyo$$
Suppose now that $\Phi$ is a free action; then the morphism
$\Phi_{y\ast}\Colon\A(G)^{\circ}\rightarrow\B(Y)^{\circ}$ is
injective which implies immediately, thanks to $\relena$ and
$\reldyo$, that $\tilde{\Phi}_{\ast}$ is injective on all 
group-like and primitive elements. By Proposition 2.17.1 of 
{\kost}, this is a necessary and sufficient condition for the 
morphism $\tilde{\Phi}_{\ast}$ to be injective on the whole 
graded coalgebra $\A(G)^{\circ}$. The converse is immediate 
again by $\relena$ and $\reldyo$.\qed

\vskip0.3cm

Consider now a right action
$\Phi\Colon(Y,\B)\times(G,\A)\rightarrow(Y,\B)$; by Theorem
{\thewrima}, we have a linear map ${\bit I}_{\Phi}\Colon\frak
g\rightarrow\Der\B(Y)$ defined as ${\bit
I}_{\Phi}(a)=(\Phi^{\ast})_{a}$. We thus obtain a subspace
$\Der_{\Phi}\B(Y)=\im{\bit I}_{\Phi}$ of the Lie superalgebra of
derivations on $\B(Y)$. As a matter of fact, $\Der_{\Phi}\B(Y)$ is a
graded Lie subalgebra of $\Der\B(Y)$. Indeed, one readily verifies
that for all $a,b\an\frak g$ we have
$[(\Phi^{\ast})_{a},(\Phi_{\ast})_{b}]=(\Phi^{\ast})_{[a,b]}$, which
means that ${\bit I}_{\Phi}([a,b])=[{\bit I}_{\Phi}(a),{\bit
I}_{\Phi}(b)]$. The following theorem provides an important property
of free actions on graded manifolds.

\math{Theorem.}{\freeactiontheorem}{\sl Let $\Phi$ be a free right action
of the graded Lie group $(G,\A)$ on the graded manifold $(Y,\B)$, 
$\dim(G,\A)=(m,n)$. Then $\Phi$ induces a
regular and involutive graded distribution $\E$ on $(Y,\B)$ of
dimension $(m,n)$.}

\undertext{\it Proof}.
 
$\bullet$ {\it Step 1}. Let us first calculate the kernel of the Lie
superalgebra morphism ${\bit I}_{\Phi}\Colon\frak
g\rightarrow\Der\B(Y)$ when $\Phi$ is a free action. To this end, the
following general property of actions is useful:
$$\widetilde{(\Phi^{\ast})_{a}}(y)=\Phi_{y\ast}(a),\;\forall
y\an Y,\;\forall a\an\frak g.\eqn\generalprop$$
For the proof of $\generalprop$, we note only that, by relation
$\derivation$, $\widetilde{(\Phi^{\ast})_{a}}(y)=
\delta_{y}\comp(\Phi^{\ast})_{a}=\Phi_{\ast}(\delta_{y}\otimes a)$.
Suppose now that ${\bit I}_{\Phi}(a)=0\Leftrightarrow(\Phi^{\ast})_{a}=0$;
by $\generalprop$, this implies that
$\Phi_{y\ast}(a)=0,\forall y\an Y$. Since $\Phi$ is a free
action, we know by Definition {\freeaction}, that $\Phi_{y\ast}$ is
injective for all $y\an Y$, which implies that $a=0$. As a result, 
$\ker{\bit I}_{\Phi}=0$, or ${\bit I}_{\Phi}$ is injective; hence, 
$\Der_{\Phi}\B(Y)$ is a graded Lie subalgebra of $\Der\B(Y)$ whose
even and odd dimensions are $m$ and $n$ respectively:
$(\Der_{\Phi}\B(Y))_{0}\cong\frak g_{0}$, 
$(\Der_{\Phi}\B(Y))_{1}\cong\frak g_{1}$.

$\bullet$ {\it Step 2}. Let now $\Der_{\Phi}\B$ be the correspondence 
$U\rightarrow\Der_{\Phi}\B(U)$, where $\Der_{\Phi}\B(U)=P_{YU}
\big(\Der_{\Phi}\B(Y)\big)$, $U\subset Y$ and 
$P_{UV}\Colon\Der\B(U)\rightarrow\Der\B(V)$ are the restriction 
maps for the sheaf $\Der\B$. Clearly, 
$\Der_{\Phi}\B$ is a subpresheaf of $\Der\B$. Consider now the
subpresheaf $\E=\B\cdot\Der_{\Phi}\B$ of $\Der\B$, $\E(U)=\B(U)\cdot
\Der_{\Phi}\B(U)=P_{YU}\big(\B(Y)\cdot\Der_{\Phi}\B(Y)\big)$. $\E(U)$
is the set of finite linear combinations of elements of $
\Der_{\Phi}\B(U)$ with coefficients in $\B(U)$. In order to prove 
that $\E$ is a sheaf, let us consider an open $U\subset Y$, an open 
covering $\{U_{\alpha}\}_{\alpha\an\Lambda}$ of $U$ and elements
$D_{\alpha}\an\E(U_{\alpha})$ such that $P_{U_{\alpha}U_{\alpha\beta}}
(D_{\alpha})=P_{U_{\beta}U_{\alpha\beta}}(D_{\beta})$,
$\forall\alpha,\beta\an\Lambda$ when $U_{\alpha\beta}=U_{\alpha}\cap
U_{\beta}\neq\kenosyn$. Then by the sheaf properties of $\Der\B$,
there exists an element $D\an\Der\B(U)$ such that
$P_{UU_{\alpha}}(D)=D_{\alpha}$. But if we write 
$D_{\alpha}=\sum_{i}f^{i}_{\alpha}P_{YU_{\alpha}}
(\Phi^{\ast})_{e_{i}}$, $f^{i}_{\alpha}\an\B(U_{\alpha})$  and
$\{e_{i}\}$ is a basis of $\frak g$, then by step 1, we find easily
that $f^{i}_{\alpha}=f^{i}|_{U_{\alpha}}$, $f^{i}\an\B(U)$, because
$D_{\alpha}$'s coincide on the intersections $U_{\alpha\beta}$. This
means that $D_{\alpha}=P_{YU_{\alpha}}(E)$, where
$E=\sum_{i}F^{i}(\Phi^{\ast})_{e_{i}}$ with $F^{i}\an\B(Y)$ such that
$F^{i}|_{U}=f^{i}$, $\forall i$. It is then immediate that
$D=P_{YU}(E)\an\E(U)$. 

$\bullet$ {\it Step 3}. It is evident that the sheaf $\E$ 
previously constructed, has the properties of a graded distribution.
In fact, $\E(U)$ is a graded submodule of $\Der\B(U)$ of dimension 
$(m,n)$ for each open $U\subset Y$. This distribution is
clearly regular thanks to relation $\generalprop$ and to the fact that
the action is free. It remains to show that
it is involutive. To this end, consider two elements
$\xi=\sum_{i}f^{i}P_{YU}(\Phi^{\ast})_{a^{i}}$ and 
$\eta=\sum_{j}g^{j}P_{YU}(\Phi^{\ast})_{b^{j}}$ of $\E(U)$,
with $f^{i},g^{j}\an\B(U), a^{i},b^{j}\an\frak g$. 
Then, direct calculation shows that 
$$\eqalign{[\xi,\eta]&=\sum_{i,j}f^{i}
\big(P_{YU}(\Phi^{\ast})_{a^{i}}
g^{j}\big)P_{YU}(\Phi^{\ast})_{b^{j}}\cr
\hfill&\hskip0.4cm-(-1)^{|\xi||\eta|}\sum_{i,j}g^{j}
\big(P_{YU}(\Phi^{\ast})_{b^{j}}f^{i}\big)
P_{YU}(\Phi^{\ast})_{a^{i}}\cr
\hfill&\hskip0.4cm+\sum_{i,j}(-1)^{|a^{i}||g^{j}|}f^{i}g^{j}
P_{YU}(\Phi^{\ast})_{[a^{i},b^{j}]},\cr}$$
from which the involutivity is evident.\qed

\vskip0.3cm

We focus now our attention on graded quotient structures defined by
equivalence relations on graded manifolds (see {\alm} for a general
treatment on this subject). A special case of equivalence relation is
provided by the action of a graded Lie group on a graded manifold and
this will be the interesting one for us. 

\math{Definition.}{\regularaction}{\sl We call a right action
$\Phi\Colon(Y,\B)\times(G,\A)\rightarrow(Y,\B)$ regular if the
morphism $\tilde{\Phi}=(\Phi\times\pi_{1})\comp\Delta\Colon(Y,\B)
\times(G,\A)\rightarrow(Y,\B)\times(Y,\B)$ defines
$(Y,\B)\times(G,\A)$ as a closed graded submanifold of
$(Y,\B)\times(Y,\B)$.}

\vskip0.3cm

Recall here, {\kost}, that $(R,\D)$ is a graded submanifold of
$(Y,\B)$ if $\D(R)^{\circ}\subset\B(Y)^{\circ}$ and there exists a
morphism of graded manifolds $i\Colon(R,\D)\rightarrow(Y,\B)$ such
that $i_{\ast}\Colon\D(R)^{\circ}\rightarrow\B(Y)^{\circ}$ is simply
the inclusion; $(R,\D)$ will be called closed if, furthermore, 
$\dim(R,\D)<\dim(Y,\B)$. Then, the action is regular if the subset
$\tilde{\Phi}_{\ast}\big(\B(Y)^{\circ}\otimes\A(G)^{\circ}\big)\subset
\B(Y)^{\circ}\otimes\B(Y)^{\circ}$ defines a graded submanifold of
$(Y,\B)\times(Y,\B)$ in the sense of Kostant, {\kost}. 
The following theorem generalizes in a natural way to the graded 
case, a fundamental result about quotients defined by actions in 
ordinary manifold theory.

\math{Theorem.}{\regularactiontheorem}{\sl  
The action $\Phi\Colon(G,\A)\times(Y,\B)\rightarrow(Y,\B)$
is regular if and only if the quotient $(Y/G,\B/\A)$ is a 
graded manifold.}

\undertext{\it Proof}. Thanks to Theorem 2.6 of {\alm},
it suffices to prove that the projections
$p_{i}\Colon(Y,\B)\times(Y,\B)\rightarrow(Y,\B), i=1,2$ on the first
and second factors restricted to the image of $(Y,\B)\times(G,\A)$ 
under $\tilde{\Phi}$ are submersions.
In other words, we must show that the morphisms of graded coalgebras
$p_{i\ast}\comp\tilde{\Phi}_{\ast}\Colon\B(Y)^{\circ}\otimes
\A(G)^{\circ}\rightarrow\B(Y)^{\circ},i=1,2$ restricted to primitive
elements are surjective. 

Consider an arbitrary primitive element $V=u\otimes a+b\otimes w$, 
for $a=\delta_{g}\an\A(G)^{\circ},b=\delta_{y}\an\B(Y)^{\circ}$
group-like and $u\an T_{y}(Y,\B)$, $w\an T_{g}(G,\A)$ primitive. 
Using relation $\reldyo$, we find easily:
$p_{1\ast}\tilde{\Phi}_{\ast}(V)=\Phi_{g\ast}(u)+\Phi_{y\ast}(w)$
and $p_{2\ast}\tilde{\Phi}_{\ast}(V)=u$, which proves that 
$p_{i}\comp\tilde{\Phi}$ are submersions, $i=1,2$.\qed

\vskip0.3cm

Note here that if $U\subset Y/G$ is an open subset, then the sheaf 
$\B/\A$ is given by 
$$(\B/\A)(U)=\{f\an\B(\check{\pi}^{-1}(U))\,|\,\Phi^{\ast}f=
f\otimes\ena_{\script A}\},\eqn\phiastf$$
where $\check{\pi}\Colon Y\rightarrow Y/G$ is the projection, {\alm}. 
Furthermore, the dimension of the quotient graded manifold 
$(Y/G,\B/\A)$ is equal to $\dim(Y/G,\B/\A)=2\dim(Y,\B)-
\dim(\im\tilde{\Phi})$, where $\im\tilde{\Phi}$
denotes the closed graded submanifold defined by $\tilde{\Phi}$. When
$\Phi$ is a free action, then by Proposition {\freeactionequiv}, we
take that $\dim(Y/G,\B/\A)=\dim(Y,\B)-\dim(G,\A)$.

We make finally some comments about graded isotropy subgroups
recalling their construction from {\kost}, but in a more concise way. 
Consider a right action $\Phi\Colon(Y,\B)\times(G,\A)
\rightarrow(Y,\B)$ and $b\an\B(Y)^{\circ}$ a group-like element, 
$b=\delta_{y}$, $y\an Y$. 
Let $H_{y}(G,\frak g)$ be the set of elements $a\an\A(G)^{\circ}$ 
with the property 
$$\Phi_{\ast}(\delta_{y}\otimes a)=\epsilon_{\script A}^{\circ}
(a)\delta_{y}.\eqn\isotropyrelationena$$
Let $H_{y}(G,\frak g)\cap G=G_{y}$ and 
$H_{y}(G,\frak g)\cap\frak g=\frak g_{y}$, then $(\frak g_{y})_{0}$ is
the Lie algebra of $G_{y}$. It is then clear that we can form the 
Lie-Hopf algebra $\R(G_{y})\crossed E(\frak g_{y})$ because $\frak
g_{y}$ is stable under the adjoint action of $G_{y}$, see Proposition 
{\adjointactionprop} and relation $\isotropyrelationena$. By Proposition 
3.8.3 of {\kost}, $\R(G_{y})\crossed E(\frak g_{y})$ corresponds to a 
graded Lie subgroup of $(G,\A)$. We denote this subgroup by $(G_{y},\A_{\kern
-1.5pt y})$ and call it graded isotropy subgroup of $(G,\A)$ at the
point $y$.

\chapter{Graded principal bundles}

Graded principal bundles were first introduced in {\alm}, {\almdyo}.
Here, we discuss this notion with slight modifications suggested by
the requirement that the definition of graded principal bundles
reproduces well the ordinary principal bundles.

\math{Definition.}{\gpb}{\sl A graded principal bundle over a graded
manifold $(X,\C)$ consists of a graded manifold $(Y,\B)$ and an action
$\Phi$ of a graded Lie group $(G,\A)$ on $(Y,\B)$ with the following
properties:
\item{1.} $\Phi$ is a free right action
\item{2.} the quotient $(Y/G,\B/\A)$ is a graded manifold, isomorphic
to $(X,\C)$, such that the natural projection $\pi\Colon(Y,\B)
\rightarrow(X,\C)$ is a submersion
\item{3.} $(Y,\B)$ is locally trivial that is, for each open $U\subset
X$, there exists an isomorphism of graded manifolds
$\phi\Colon(V,\B|_{V})\rightarrow(U\times G,\C|_{U}\opi\A)$, 
$V=\pi_{\ast}^{-1}(U)\subset Y$, such that the isomorphism
$\phi^{\ast}$ of graded algebras is a morphism of $\A(G)$-comodules,
where the $\A(G)$-comodule structures on $\C(U)\opi\A(G)$ and $\B(V)$
are given by $id\otimes\da$ and $\Phi^{\ast}$ respectively.
Furthermore, we require that $\phi^{\ast}=m_{\script
B}\comp(\pi^{\ast}\otimes\psi^{\ast})$, where
$\psi\Colon(V,\B|_{V})\rightarrow(G,\A)$ is a morphism of graded
manifolds.

}

\vskip0.3cm

The fact that $\phi^{\ast}$ is a morphism of $\A(G)$-comodules,
that is, $$(\phi^{\ast}\otimes id)\comp(id\otimes\da)=\Phi^{\ast}
\comp\phi^{\ast}\eqn\comodmorphi$$
implies that $\psi^{\ast}$ is also a morphism of $\A(G)$-comodules:
$$(\psi^{\ast}\otimes id)\comp\da=\Phi^{\ast}\comp\psi^{\ast}.
\eqn\comodmorpsi$$
One easily verifies that the underlying differentiable manifolds of
Definition {\gpb} form an ordinary principal bundle, and further, if 
the graded manifolds become trivial, in the sense that 
$\A=C^{\infty}_{G},\B=C^{\infty}_{Y}, \C=C^{\infty}_{X}$, 
then we obtain the definition of an ordinary principal bundle. 
We refer the reader to {\kobay} for a general and systematic treatment
on the subject of principal bundles in differential geometry. 

Let us now compute the graded isotropy subgroups in the case 
where the action $\Phi$ is free (for example, this is the case 
of the graded principal bundle).
The Lie group $G_{y}$ which is defined as $G_{y}=\{g\an G\,|\,\Phi_{\ast}
(b\otimes\delta_{g})=b\}$ is equal to $e$ because the action of $G$ on 
$Y$ is free. On the other hand, $\frak g_{y}=\{a\an\frak g\,|\,
\Phi_{\ast}(b\otimes a)=0\}=0$, again because the action is 
free (the morphism $\Phi_{y\ast}$ is injective). Hence, in this case 
the graded isotropy subgroup $(G_{y},\A_{\kern-1.5pt y})$ is simply 
$(e,{\script R})$, where $\script R$ is the trivial sheaf over the 
identity $e\an G$, ${\script R}(e)=\R$.

In order to calculate the quotient graded manifold
$(G/G_{y},\A/\A_{\kern-1.5pt y})$ which represents the orbit of the
point $y$ under the action of $(G,\A)$, we need the expression of 
the canonical right action $\Phi$ of the subgroup $(G_{y},\A_{\kern-1.5pt
y})$ on $(G,\A)$. If $i\Colon(G_{y},\A_{\kern-1.5pt
y})\rightarrow(G,\A)$ is the inclusion, we have:
$\Phi^{\ast}=(id\otimes i^{\ast})\comp\da$ and for the case of the free 
action, where $(G_{y},\A_{\kern-1.5pt y})=(e,{\script R})$, one finds that
$i^{\ast}=\epsilon_{\script A}$ and finally $\Phi^{\ast}=id$. In view of
relation $\phiastf$, it is straightforward that $\A/\A_{\kern-1.5pt y}=\A$.
Therefore,

\math{Property.}{\idiotita}{\sl The orbits of a free action of
$(G,\A)$ are always isomorphic as graded manifolds to $(G,\A)$.}

\vskip0.3cm

For the case now of the graded principal bundle, the orbit
$(\O_{y},\B_{y})$ of $(G,\A)$
through $y\an Y$ will be called fibre of $(Y,\B)$ over
$x=\pi_{\ast}|_{Y}(y)\an X$. Using the graded version of the submersion
theorem, {\ruipma}, one can justify this terminology as follows. 
The pre-image $\pi^{-1}(x,{\script R})$ of the closed graded submanifold
$(x,{\script R})\hookrightarrow(X,\C)$ is a closed graded submanifold of
$(Y,\B)$ whose underlying differentiable manifold is
$\pi_{\ast}^{-1}(x)$. So, if we write $\pi^{-1}(x,{\script R})=
(\pi_{\ast}^{-1}(x),\D)$, then for each $z\an\pi_{\ast}^{-1}(x)$ we have: 
$T_{z}(\pi_{\ast}^{-1}(x),\D)=\ker T_{z}\pi$.  We know already that
$\pi_{\ast}^{-1}(x)=\O_{y}$, the orbit under $G$ of a point $y\an Y$
such that $\pi_{\ast}|_{Y}(y)=x$. Furthermore, if
$\delta_{z}=\Phi_{\ast}(\delta_{y}\otimes\delta_{g})$, $g\an G$ and
$v\an T_{g}(G,\A)$, then $V=\Phi_{\ast}(\delta_{y}\otimes v)\an
T_{z}(\O_{y},\B_{y})$ and $\pi_{\ast}(V)=0$. Consequently, 
$T_{z}(\O_{y},\B_{y})\subset\ker T_{z}\pi$. By a  simple argument
on dimensions, we obtain that $T_{z}(\O_{y},\B_{y})=
T_{z}(\pi_{\ast}^{-1}(x),\D)$. We conclude that the tangent bundles 
of $\pi^{-1}(x,{\script R})$ and $(\O_{y},\B_{y})$ are identical; but
then, Theorem 2.16 of {\kost} tells us that these graded manifolds 
coincide. 

Next, we discuss an elementary example of graded principal bundle, the
product bundle. In this case, one can directly verify the axioms of 
Definition {\gpb}. Nevertheless, there exist also graded
principal bundles for which Definition {\gpb} cannot be directly
applied, even though this is possible for the corresponding ordinary
principal bundles. For such cases, one may use an equivalent
definition of the graded principal bundle, see next section. 

\math{Example.}{\paradeigmaena}{Consider a graded manifold $(X,\C)$, 
a graded Lie group $(G,\A)$ and their product
$(Y,\B)=(X,\C)\times(G,\A)$. One has a canonical right action 
$\Phi\Colon(Y,\B)\times(G,\A)\rightarrow(Y,\B)$ defined as 
$\Phi^{\ast}=id\otimes\da$. This action is free: if $\delta_{y}=
\delta_{x}\otimes\delta_{g}\an\C(X)^{\circ}\otimes\A(G)^{\circ}$ is 
group-like and $a\an\A(G)^{\circ}$, then $\Phi_{y\ast}(a)=\Phi_{\ast}
(\delta_{x}\otimes\delta_{g}\otimes a)=
\delta_{x}\otimes(\delta_{g}\odot a)=\delta_{x}\otimes L_{g\ast}(a)$, 
which implies that $\Phi_{y\ast}$ is an injective
morphism of graded coalgebras, because $L_{g\ast}$ is an 
isomorphism. Evidently, the quotient $Y/G$ is equal to $X$ and 
the sheaf $\B/\A$ 
over $X$ is given by the elements $f\an\C(U)\opi\A(G)$ for which 
$\Phi^{\ast}f=f\otimes\ena_{\script A}$, $U\subset X$ open. 
If we decompose $f$ as $f=\sum_{i}f_{i}\otimes h_{i}$, $f_{i}\an\C(U), 
h_{i}\an\A(G)$, we take easily $\da h_{i}=h_{i}\otimes
\ena_{\script A}$, hence $h_{i}$'s are such that $\epsilon_{\script
A}(h_{i})\ena_{\script A}=h_{i}$. We conclude that $f$ is of the form 
$f=\sum_{i}\epsilon_{\script A}(h_{i})f_{i}\otimes\ena_{\script
A}=f_{\script C}\otimes\ena_{\script A}$ and finally 
$(\B/\A)(U)\cong\C(U)$, which proves that the quotient $(Y/G,\B/\A)$ 
is isomorphic to $(X,\C)$. Further, the identity map
$\phi^{\ast}=id\Colon\C(U)\opi\A(G)\rightarrow\C(U)\opi\A(G)$ 
admits the decomposition of Definition {\gpb} (it suffices to choose
$\pi=\pi_{1},\psi=\pi_{2}$, the projections on the first and second
factors respectively) and satisfies trivially the relation
$\comodmorphi$.}

\chapter{The geometry of graded principal bundles}

In this section we analyze three aspects of the geometry of graded
principal bundles: the relation between the sheaf of vertical
derivations and the graded distribution induced by the action of the
structure group, a criterion of global triviality of the graded
principal bundle, and, finally, a way to reformulate Definition {\gpb}
avoiding the use of local trivializations. 

For this and the subsequent sections, we will adopt the following
notation in order to simplify the discussion: if $a\an\frak g$ and
$V\subset Y$ is an open, then the restriction 
$P_{YV}(\Phi^{\ast})_{a}$ (Theorem {\freeactiontheorem}) will be 
simply denoted by $(\Phi^{\ast})_{a}$.

It is well-known that if $Y(X,G)$ is an ordinary principal bundle, 
then the set of vertical vectors at $y\an Y$ is equal to the set 
of induced vectors at the same point. In the previous section, 
we saw that the same is true for graded principal bundles; 
however, it is not evident that this property remains valid for 
the sheaves of vertical and induced derivations. Nonetheless, 
as the following theorem confirms, this is indeed the case. 

\math{Theorem.}{\ekerpi}{\sl Let $(Y,\B)$ be a graded principal bundle
over $(X,\C)$ with structure group $(G,\A)$. If $\E$ is the natural
graded distribution induced by the free action of $(G,\A)$ on
$(Y,\B)$, then $\E$ is equal to the sheaf of vertical derivations,
$\E=\Ver(\pi_{\ast},\B)$.}

\undertext{\it Proof}. We show first that
$\E\subset\Ver(\pi_{\ast},\B)$; to this end, it is sufficient to prove
that for each $a\an\frak g$, we have $\pi_{\ast}(\Phi^{\ast})_{a}=0$. 
Indeed, if $f$ is a homogeneous element of $\C(U)$, $U\subset X$ open,
we take:
$$\pi^{\ast}\big[\pi_{\ast}(\Phi^{\ast})_{a}(f)\big]=(\Phi^{\ast})_{a}
(\pi^{\ast}f)=(id\otimes a)(\pi^{\ast}f\otimes\ena_{\script A})=0,$$ 
since $a$ is primitive with respect to $\delta_{e}$. 
Now the following argument on dimensions completes the proof. A
derivation $D\an\Ver(\pi_{\ast},\B)(V)$, $V=\pi_{\ast}^{-1}(U)$,
is characterized by the property: 
$\pi^{\ast}[(\pi_{\ast}D)f]=D(\pi^{\ast}f)=0$,
$\forall f\an\C(U)$. This means in terms of coordinates that $D$ does
not depend on the graded coordinates on $V$ obtained by pulling-back
the graded coordinates of $U$ via $\pi^{\ast}$. Using the fact that
$\pi^{\ast}$ is an injection, we find that the dimension of 
$\Ver(\pi_{\ast},\B)(V)$ equals to $\dim(Y,\B)-\dim(X,\C)=\dim(G,\A)=
\dim\E(V)$.\qed

\vskip0.3cm

The fact that a local trivialization $\phi$ is an isomorphism
of $\A(G)$-comodules is expressed by relation $\comodmorphi$ but it is
also reflected in the induced derivations. The following lemma makes
this precise, providing a relation between them.

\math{Lemma.}{\indvectorfields}{\sl If $\phi\Colon(V,\B|_{V})
\rightarrow(U,\C|_{U})\times(G,\A)$, $V=\pi_{\ast}^{-1}(U)$, 
is a local trivialization of $(Y,\B)$,
then the following relation is true for each $a\an\frak g$:
$$\phi_{\ast}(\Phi^{\ast})_{a}=id\otimes (R^{\ast})_{a}.$$}
\indent\undertext{\it Proof}. We show first that $\phi_{\ast}
(\Phi^{\ast})_{a}\an\Der\A(G)$. Indeed, if $f_{\script C}\an\C(U)$, we
take:
$$\phi_{\ast}(\Phi^{\ast})_{a}(f_{\script C}\otimes\ena_{\script
A})=\big((\phi^{-1})^{\ast}\comp(\Phi^{\ast})_{a}\comp\phi^{\ast}
\big)(f_{\script C}\otimes\ena_{\script A})=
(\phi^{\ast})^{-1}(\Phi^{\ast})_{a}(\pi^{\ast}f_{\script C})=0.$$
It is then sufficient to calculate $\phi_{\ast}(\Phi^{\ast})_{a}$ 
on elements of $\C(U)\opi\A(G)$ of the form $\ena_{\script C}
\otimes f_{\script A}$ for $f_{\script A}\an\A(G)$. Taking into
account relation $\comodmorpsi$ and if  $\da
f_{\script A}=\sum_{i}I^{i}f_{\script A}\otimes J^{i}f_{\script A}$,
one finds:
$$\eqalign{\phi_{\ast}(\Phi^{\ast})_{a}(\ena_{\script C}\otimes
f_{\script A})&=(\phi^{-1})^{\ast}(id\otimes a)(\psi^{\ast}\otimes id)
\da f_{\script A}\cr
\hfill&=\sum_{i}(-1)^{|a||I^{i}f_{\script A}|}
(\phi^{\ast})^{-1}\psi^{\ast}(I^{i}f_{\script A})a(J^{i}f_{\script
A})\cr
\hfill&=\sum_{i}(-1)^{|a||I^{i}f_{\script A}|}(\ena_{\script C}\otimes
I^{i}f_{\script A})a(J^{i}f_{\script A})\cr
\hfill&=\ena_{\script C}\otimes(id\otimes a)\da f_{\script A}
=(id\otimes(R^{\ast})_{a})(\ena_{\script C}\otimes f_{\script A}),\cr}$$
where we have used that $\phi^{\ast}(\ena_{\script C}\otimes
I^{i}f_{\script A})=\psi^{\ast}(I^{i}f_{\script A})$ implies
$(\phi^{-1})^{\ast}\psi^{\ast}(I^{i}f_{\script A})=
\ena_{\script C}\otimes I^{i}f_{\script A}$.\qed

\vskip0.3cm

Next we discuss the notion of section of a graded principal bundle
and we show that graded and ordinary sections exhibit several
analogous properties.

\math{Definition.}{\gradedsectiondef}{\sl Let $U\subset X$ be an open
on the base manifold $(X,\C)$ of a graded principal bundle $(Y,\B)$. 
We call graded section of $(Y,\B)$ on $U$ a morphism of graded
manifolds $s\Colon(U,\C|_{U})\rightarrow(Y,\B)$ having the property
$s^{\ast}\comp\pi^{\ast}=id$. We write also $\pi\comp
s=id\Colon(U,\C|_{U})\rightarrow(U,\C|_{U})$.}

\vskip0.3cm

A  first property of graded sections is that to each local
trivialization, one can associate in a canonical way a graded section. 
More precisely:

\math{Lemma.}{\gradedsectionlem}{\sl Let $\phi\Colon(V,\B|_{V})
\rightarrow(U,\C|_{U})\times(G,\A)$, $V=\pi_{\ast}^{-1}(U)$, be a local
trivialization. Then, if $E\Colon\C(U)\opi\A(G)\rightarrow\C(U)$ is
defined as $E(f_{\script C}\otimes f_{\script A})=\delta_{e}
(f_{\script A})f_{\script C}$, the map $E\comp(\phi^{-1})^{\ast}
\Colon\B(V)\rightarrow\C(U)$ defines a morphism of graded manifolds 
with the properties of a graded section.}

\undertext{\it Proof}. The fact that $E$ is a morphism of graded
manifolds is evident because we may write $E=id\otimes\delta_{e}$.
Therefore, there exists a morphism of graded manifolds
$s\Colon(U,\C|_{U})\rightarrow(Y,\B)$ such that
$s^{\ast}=E\comp(\phi^{-1})^{\ast}$. Now if $f_{\script C}\an\C(U)$,
we have: $(s^{\ast}\comp\pi^{\ast})(f_{\script
C})=E\big((\phi^{-1})^{\ast}\pi^{\ast}f_{\script C}\big)$ and since
$\phi^{\ast}(f_{\script C}\otimes\ena_{\script
A})=\pi^{\ast}f_{\script C}$, we finally obtain $s^{\ast}\comp
\pi^{\ast}=id$.\qed

\vskip0.3cm 

Conversely now, consider a graded section $s\Colon(U,\C|_{U})
\rightarrow(Y,\B)$. We wish to show that there exists a local 
trivialization $\phi\Colon(V,\B|_{V})\rightarrow(U,\C|_{U})
\times(G,\A)$ canonically associated to $s$,
$V=\pi_{\ast}^{-1}(U)\subset Y$. To this end, we first define
a morphism of graded algebras $\tilde{\phi}^{\ast}\Colon\B(V)
\rightarrow\C(U)\opi\A(G)$ by $\tilde{\phi}^{\ast}=(s^{\ast}
\otimes id)\comp\Phi^{\ast}$. Let $\tilde{\phi}_{\ast}=\Phi_{\ast}
\comp(s_{\ast}\otimes id)\Colon\C(U)^{\circ}\otimes\A(G)^{\circ}
\rightarrow\B(V)^{\circ}$ be the corresponding morphism of graded 
coalgebras. Clearly, the differentiable mapping 
$\tilde{\phi}_{\ast}|_{U\times G}\Colon U\times G\rightarrow V$ is 
bijective. Consider now the tangent of $\phi$ at the arbitrary point 
$(x,g)\an U\times G$. If $z=u\otimes\delta_{g}+\delta_{x}\otimes w
\an T_{(x,g)}(U\times G,\C|_{U}\opi\A)$ is a tangent vector at
$(x,g)$, $u\an T_{x}(U,\C|_{U})$, $w=T_{g}(G,\A)$, then:
 $$\tilde{\phi}_{\ast}(z)=\Phi_{g\ast}(s_{\ast}u)+
(\Phi_{s_{\ast}x})_{\ast}(w).$$
This implies that $T_{(x,g)}\tilde{\phi}$ is injective, because
$\Phi_{g\ast}$ is an isomorphism and $\Phi_{y\ast}$ is injective for
each $y\an Y$ (the action is free). Thus, $T_{(x,g)}\tilde{\phi}$
is an injection between two vector spaces of the same dimension and
hence an isomorphism. Using now Theorem 2.16 of {\kost}, we conclude
that $\tilde{\phi}$ is an isomorphism of graded manifolds.

Let now $\pi_{1}\Colon(U,\C|_{U})\times(G,\A)\rightarrow(U,\C|_{U})$ and 
$\pi_{2}\Colon(U,\C|_{U})\times(G,\A)\rightarrow(G,\A)$ be the
projections. Then:

\math{Proposition.}{\localtr}{\sl Let $\phi\Colon(V,\B|_{V})
\rightarrow(U,\C|_{U})\times(G,\A)$, $V=\pi_{\ast}^{-1}(U)$, 
be the morphism of graded manifolds defined as $\phi^{\ast}=
m_{\script B}\comp(\pi^{\ast}\otimes\psi^{\ast})$, where $\psi^{\ast}=
(\tilde{\phi}^{\ast})^{-1}\comp\pi_{2}^{\ast}$. Then, $\phi$ is a local
trivialization of the graded principal bundle $(Y,\B)$.}

\undertext{\it Proof}. We show first that $\phi^{\ast}$ is an 
isomorphism. To this end, consider the composition 
$\tilde{\phi}^{\ast}\comp\phi^{\ast}$:
$$\eqalign{\tilde{\phi}^{\ast}\comp\phi^{\ast}&=
(s^{\ast}\otimes id)\comp\Phi^{\ast}\comp m_{\script
B}\comp(\pi^{\ast}\otimes\psi^{\ast})\cr
\hfill&=m_{\script C\script A}\comp[(s^{\ast}\otimes
id)\otimes(s^{\ast}\otimes
id)]\comp(\Phi^{\ast}\otimes\Phi^{\ast})\comp(\pi^{\ast}
\otimes\psi^{\ast})\cr
\hfill&=m_{\script C\script A}\comp[(s^{\ast}\otimes
id)\comp\Phi^{\ast}\comp\pi^{\ast}\otimes(s^{\ast}\otimes
id)\comp\Phi^{\ast}\comp\psi^{\ast}]\cr
\hfill&=m_{\script C\script A}\comp[(s^{\ast}\otimes
id)\comp\Phi^{\ast}\comp\pi^{\ast}\otimes\pi^{\ast}_{2}].\cr}$$
Using now the fact that $(s^{\ast}\otimes id)\comp\Phi^{\ast}\comp
\pi^{\ast}=\pi_{1}^{\ast}$, we find that $\tilde{\phi}^{\ast}\comp
\phi^{\ast}=id$, which proves that $\phi^{\ast}$ is also an isomorphism. 
It remains to show relation $\comodmorpsi$, or equivalently,
$\psi_{\ast}\comp\Phi_{\ast}=m_{\script A}^{\circ}\comp
(\psi_{\ast}\otimes id)$. Since $\tilde{\phi}_{\ast}$ is an isomorphism
between $\B(V)^{\circ}$ and $\C(U)^{\circ}\otimes\A(G)^{\circ}$, we
may write each element $b\an\B(V)^{\circ}$ as $b=\Phi_{\ast}
(\sum_{i}s_{\ast}c^{i}\otimes a^{i}_{0})$, for $c^{i}\an\C(U)^{\circ}$ and
$a^{i}_{0}\an\A(G)^{\circ}$. If $a\an\A(G)^{\circ}$, we take:
$$\eqalign{\psi_{\ast}\Phi_{\ast}(b\otimes a)&=
\psi_{\ast}\Phi_{\ast}\big(\Phi_{\ast}(\sum_{i}s_{\ast}c^{i}\otimes
a^{i}_{0})\otimes a\big)=\psi_{\ast}\Phi_{\ast}
\big(\sum_{i}s_{\ast}c^{i}\otimes(a^{i}_{0}\odot a)\big)\cr
\hfill&=\psi_{\ast}\Phi_{\ast}(s_{\ast}\otimes id)\big(\sum_{i}c^{i}
\otimes(a^{i}_{0}\odot a)\big)=
\pi_{2\ast}\big(\sum_{i}c^{i}\otimes(a^{i}_{0}\odot a)\big)\cr
\hfill&=\sum_{i}c^{i}(\ena_{\script C})a^{i}_{0}\odot a=
\psi_{\ast}b\odot a,\cr}$$
since $\psi_{\ast}b=\pi_{2\ast}(\sum_{i}c^{i}\otimes a^{i}_{0})=
\sum_{i}c^{i}(\ena_{\script C})a^{i}_{0}$.\qed

\vskip0.3cm

By Lemma {\gradedsectionlem} and Proposition {\localtr}, if we set
$U=X$, it is straightforward that the condition of global triviality 
of a principal bundle remains valid in the graded setting.

\math{Corollary-Theorem.}{\porismatheorima}{\sl A graded principal
bundle $(Y,\B)$ is globally isomorphic to the product $(X,\C)
\times(G,\A)$ if and only if it admits a global section
$s\Colon(X,\C)\rightarrow(Y,\B)$.}

\vskip0.3cm

We observe here that Lemma {\gradedsectionlem} and 
Proposition {\localtr} remain valid if we replace the graded principal
bundle $(Y,\B)$ by a graded manifold $(Y,\B)$ on which the graded Lie
group $(G,\A)$ acts freely to the right in such a way that the
quotient $(X,\C)=(Y/G,\B/\A)$ is a graded manifold, and the projection 
$\pi\Colon(Y,\B)\rightarrow(X,\C)$ is a submersion. In that case,
one can construct via Lemma {\gradedsectionlem} and Proposition {\localtr}
the local trivializations of Definition {\gpb}. In other words:

\math{Theorem.}{\orismoskaithewrima}{\sl A graded principal
bundle is a graded manifold $(Y,\B)$ together with a free right
action $\Phi\Colon(Y,\B)\times(G,\A)\rightarrow(Y,\B)$ of a graded Lie
group $(G,\A)$ such that: 
\item{1.} the quotient $(X,\C)=(Y/G,\B/\A)$ is a graded manifold 
\item{2.} the projection $\pi\Colon(Y,\B)\rightarrow(X,\C)$ is a
submersion. 

}

\vskip0.3cm

As an immediate application, we examine if the principal bundles
formed by Lie groups and closed Lie subgroups possess graded
analogs. 

\math{Example.}{\paradeigmadyo}{ Consider a graded Lie group $(G,\A)$ 
and a closed graded Lie subgroup $(H,\D)$ of $(G,\A)$. The natural
right action $\Phi\Colon(G,\A)\times(H,\D)\rightarrow(G,\A)$ is
given by $\Phi^{\ast}=(id\otimes i^{\ast})\comp\da$, where $i\Colon
(H,\D)\rightarrow(G,\A)$ is the inclusion. Furthermore, we know 
that the quotient $(G/H,\A/\D)$ is a graded manifold and the
projection $(G,\A)\rightarrow(G/H,\A/\D)$ is a submersion, {\kost}.
Now if $g\an G$, the morphism $\Phi_{g\ast}\Colon\D(H)^{\circ}\rightarrow
\A(G)^{\circ}$ is given by $\Phi_{g\ast}(d)=\delta_{g}\odot d=
L_{g\ast}(d)$, $\forall d\an\D(H)^{\circ}$. As a result,
$\Phi_{g\ast}$ is injective, so the action $\Phi$ is free and
Theorem {\orismoskaithewrima} holds: $(G,\A)$ is a graded
principal bundle over $(G/H,\A/\D)$ with typical fibre $(H,\D)$.}

\chapter{Lie superalgebra-valued graded differential forms} 

Let $(Y,\B)$ be a graded manifold and $\frak g$ a Lie superalgebra.
We call $\frak g$-valued graded differential form on $(Y,\B)$, an
element of $\Omega(Y,\B)\otimes\frak g$. It is clear that the set
$\OYBG=\Omega(Y,\B)\otimes\frak g$ of these forms constitutes a
$(\Z\oplus\Z_{2}\oplus\Z_{2})$-graded vector space; however, it is more
convenient to introduce a $(\Z\oplus\Z_{2})$-grading as follows:
if $\alpha\an\OYBG$ and $\deg(\alpha)=(i_{\alpha},j_{\alpha},
k_{\alpha})$ is its $(\Z\oplus\Z_{2}\oplus\Z_{2})$-degree, then we 
set $|\alpha|=(i_{\alpha},j_{\alpha}+k_{\alpha})\an\Z\oplus\Z_{2}$.
If $\{e_{i}\}$ is a basis of the Lie superalgebra $\frak g$ and
$\alpha,\beta\an\OYBG$, then we may write $\alpha=\sum_{i}\alpha^{i}
\otimes e_{i}$ and $\beta=\sum_{i}\beta^{i}\otimes e_{i}$, for
$\alpha^{i}$, $\beta^{i}\an\Omega(Y,\B)$. In the case where $\alpha$
and $\beta$ are homogeneous with $\deg(\alpha)=(i_{\alpha},j_{\alpha},
k_{\alpha})$, $\deg(\beta)=(i_{\beta},j_{\beta},k_{\beta})$, we define
the $\frak g$-valued graded differential form $[\alpha,\beta]$ of
degree $\deg([\alpha,\beta])=(i_{\alpha}+i_{\beta},j_{\alpha}+
j_{\beta},k_{\alpha}+k_{\beta})$ as:
$$[\alpha,\beta]=\sum_{i,j}(-1)^{j_{\beta}k_{\alpha}}\alpha^{i}\beta^{j}
\otimes[e_{i}^{\alpha},e_{j}^{\beta}],\eqn\formediff$$
if $\alpha=\sum_{i}\alpha^{i}\otimes e_{i}^{\alpha}$ with
$|\alpha^{i}|=(i_{\alpha},j_{\alpha})$ and
$e_{i}^{\alpha}=k_{\alpha}$; similarly for $\beta$.
We extend this definition to non-homogeneous elements by linearity.
Clearly, equation $\formediff$ gives the same result for every basis of
the Lie superalgebra $\frak g$. Thus, we have a bilinear map 
$[-,-]\Colon\Omega^{i}(Y,\B)_{j}
\otimes\frak g_{k}\times\Omega^{i\tonoss}(Y,\B)_{j\tonoss}\otimes
\frak g_{k\tonoss}\rightarrow\Omega^{i+i\tonoss}(Y,\B)_{j+j\tonoss}
\otimes\frak g_{k+k\tonoss}$ with the following properties:

\math{Proposition.}{\propert}{\sl If $\alpha,\beta,\gamma\an\OYBG$ are
homogeneous, then we have:
\vskip0.1cm
\item{1.} $[\alpha,\beta]=-(-1)^{|\alpha||\beta|}[\beta,\alpha]$
\item{2.}
$\frak S(-1)^{|\alpha||\gamma|}\big[[\alpha,\beta],\gamma\big]=0$.

In the previous relations, we have set $|\alpha||\beta|=i_{\alpha}
i_{\beta}+(j_{\alpha}+k_{\alpha})(j_{\beta}+k_{\beta})$ and $\frak S$ 
means the cyclic sum on the argument which follows.} 

\undertext{\it Proof}. Routine calculations using relation
$\formediff$.\qed

\vskip0.3cm

We realize that the space $\OYBG$ possesses the structure of a
$(\Z\oplus\Z_{2})$-graded Lie algebra, inherited from the the Lie
superalgebra structure of $\frak g$. 

The action now of elements of $\OYBG$ on derivations can be seen as
follows: if $\alpha=\sum_{i}\alpha^{i}\otimes e_{i}\an\OYBG$ is an 
$r$-form and $\xi_{1},\ldots,\xi_{r}\an\Der\B(Y)$, we set $(\xi_{1},
\ldots,\xi_{r}|\alpha)=\sum_{i}(\xi_{1},\ldots,\xi_{r}|\alpha^{i})
\otimes e_{i}$. Accordingly, one can extend the exterior differential 
$d$ to a differential on $\frak g$-valued graded differential forms,
also noted by $d$ in the following manner: if $\alpha=\sum_{i}
\alpha^{i}\otimes e_{i}$, then we set $d\alpha=\sum_{i}d\alpha^{i}
\otimes e_{i}$. The exterior differential $d\Colon\OYBG\rightarrow\OYBG$ 
defined previously, is a derivation of degree $|d|=(1,0)$. By
straightforward verification, we find that, if $\alpha$, $\beta$ are
$\frak g$-valued graded differential forms and $\alpha$ is
homogeneous, then $d[\alpha,\beta]=[d\alpha,\beta]+(-1)^{|\alpha||d|}
[\alpha,d\beta]$. 
 
In a similar way, one can extend the pull-back of graded differential forms 
under a morphism of graded manifolds $\sigma\Colon(Y,\B)\rightarrow
(Z,{\script Y})$ to a linear map $\sigma^{\ast}\Colon\Omega(Z,{\script
Y},\frak g)\rightarrow\OYBG$ which commutes with the exterior
differential, that is, $d\comp\sigma^{\ast}=\sigma^{\ast}\comp d$,
and preserves the bracket $[-,-]$: 
$\sigma^{\ast}[\alpha,\beta]=
[\sigma^{\ast}\alpha,\sigma^{\ast}\beta]$. We have analogous
generalizations for the Lie derivative. The following properties
of the bracket and the Lie derivative on $\frak g$-valued graded
differential forms will be useful; the proof proceeds by a
straightforward calculation with graded differential forms and Lie 
superalgebra elements, and is left as an exercise for the reader.

\math{Proposition.}{\idiothtalie}{\sl If $\alpha,\beta\an\OYBG$ and
$\xi\an\Der\B(Y)$ are homogeneous, then
\item{1.} $\Lie_{\xi}[\alpha,\beta]=[\Lie_{\xi}\alpha,\beta]+(-1)^{|{\bit
L}_{\xi}||\alpha|}[\alpha,\Lie_{\xi}\beta]$
\item{2.} $\big(id\otimes ad(v)\big)[\alpha,\beta]=\big[(id\otimes
ad(v))\alpha,\beta\big]+(-1)^{|v||\alpha|}
\big[\alpha,(id\otimes ad(v))\beta\big]$, $\forall v\an\frak g$.}

\vskip0.3cm

Finally, one can restate Proposition {\pullback} for the case of 
$\frak g$-valued graded differential forms. For example, if $\sigma$
is an isomorphism of graded manifolds, relation $\pullba$ becomes: 
$$(\xi_{1},\ldots,\xi_{r}|\sigma^{\ast}\alpha)
=(\sigma^{\ast}\otimes id)(\sigma_{\ast}\xi_{1},\ldots,\sigma_{\ast}
\xi_{r}|\alpha).\eqn\neostypos$$
One can also define multilinear forms on the tangent spaces of a
graded manifold taking its values in the Lie superalgebra $\frak g$,
using a simple modification of $\multiliforms$.

\chapter{Graded connections}

We have seen that on each graded principal bundle there always 
exists a natural distribution induced by the action of the structure
group which is equal to the sheaf of vertical derivations. The choice of a
connection is essentially the choice of a complementary distribution.
More precisely:

\math{Definition.}{\connexiong}{\sl Let $(Y,\B)$ be a graded principal
bundle with structure group $(G,\A)$, over the graded manifold
$(X,\C)$. A graded connection on $(Y,\B)$ is a regular distribution
$\H\subset\Der\B$ of dimension $\dim\H=\dim(X,\C)$ such that:
\item{1.} $\H\oplus\Ver(\pi_{\ast},\B)=\Der\B$,
$\pi\Colon(Y,\B)\rightarrow(X,\C)$ is the projection
\item{2.} $\H$ is $(G,\A)$-invariant.

}

\vskip0.3cm
 
Let us explain the second statement in the previous definition: 
$\H$ will be called $(G,\A)$-invariant if, for each open 
$U\subset X$ and $D\an\H(\pi_{\ast}^{-1}(U))$, the derivations
$\Phi_{g}^{\ast}D$ and $[(\Phi^{\ast})_{a},D]=\Lie_{(\Phi^{\ast})_{a}}D$
belong also to $\H(\pi_{\ast}^{-1}(U))$, $\forall g\an G$,
$\forall a\an\frak g$. In order to put these condition in a more
compact form, we introduce the following notation
$$(\Phi^{\ast})_{a}D=\left\{\matrix{\Phi^{\ast}_{g}D,\quad\hbox{if}\quad
a=\delta_{g}, g\an G\cr
\Lie_{(\Phi^{\ast})_{a}}D,\quad\hbox{if}\quad a\an\frak g.\hfill\cr}\right.$$
Then, $\H$ will be $(G,\A)$-invariant if 
$(\Phi^{\ast})_{a}D\an\H(\pi_{\ast}^{-1}(U))$, for each element $a$ 
group-like or primitive with respect to $\delta_{e}$. We can now
reformulate this notion in terms of $\frak g$-valued graded
differential forms. 

Given the graded connection $\H$ on $(Y,\B)$, we have:
$\H(Y)\oplus\E(Y)=\Der\B(Y)$, $\E=\Ver(\pi_{\ast},\B)$ 
(see Theorem {\ekerpi}), and each derivation $\xi\an\Der\B(Y)$ 
decomposes as $\xi=\xi^{H}+\sum_{i}f^{i}(\Phi^{\ast})_{e_{i}}$, 
where $\xi^{H}\an\H(Y)$, $\{e_{i}\}$ is a basis of $\frak g$ 
and $f^{i}\an\B(Y)$. Then, we define
a 1-form $\biomega\an\Omega^{1}(Y,\B,\frak g)$ setting
$(\xi|\biomega)=\sum_{i}f^{i}\otimes e_{i}$. Let us now calculate
$(\Phi^{\ast})_{a}\biomega$, where by definition 
$$(\Phi^{\ast})_{a}\biomega=\left\{\matrix{\Phi^{\ast}_{g}\biomega,
\quad\hbox{if}\quad a=\delta_{g}, g\an G\cr
\Lie_{(\Phi^{\ast})_{a}}\biomega,\quad\hbox{if}\quad a\an\frak g.\hfill\cr}
\right.$$
Consider first the case $a=\delta_{g}$. Then, by equation $\neostypos$
we take
$$(\xi|\Phi_{g}^{\ast}\biomega)=
(\Phi_{g}^{\ast}\otimes id)(\Phi_{g\ast}\big(\sum_{i}f^{i}
(\Phi^{\ast})_{e_{i}}\big)|\biomega)$$
and using the fact that $\Phi_{g\ast}(\Phi^{\ast})_{e_{i}}=
(\Phi^{\ast})_{AD_{g^{-1}\ast}(e_{i})}$, we obtain:
$$(\xi|\Phi_{g}^{\ast}\biomega)=(\Phi_{g}^{\ast}\otimes id)
\big(\sum_{i}\Phi^{\ast}_{g^{-1}}f^{i}\otimes
AD_{g^{-1}\ast}(e_{i})\big)=(id\otimes AD_{g^{-1}\ast})(\xi|\biomega).
$$ or
$$\Phi_{g}^{\ast}\biomega=(id\otimes AD_{g^{-1}\ast})\comp\biomega.
\eqn\graconnexion$$
If now $a\an\frak g$, we have:
$$(\xi|\Lie_{(\Phi^{\ast})_{a}}\biomega)=(-1)^{|a||\xi|}
\Big[\big((\Phi^{\ast})_{a}\otimes id\big)(\xi|\biomega)-
\big([(\Phi^{\ast})_{a},\xi]|\biomega\big)\Big]\eqn\xrhsimo$$
and it is sufficient to examine two cases:
\item{1.} $\xi=\xi^{H}$ (horizontal derivation):
it is then immediate that $(\xi|\Lie_{(\Phi^{\ast})_{a}}\biomega)=0$ 
\item{2.} $\xi=(\Phi^{\ast})_{b}, b\an\frak g$ 
(vertical derivation): $$(\xi|\Lie_{(\Phi^{\ast})_{a}}\biomega)=
(-1)^{|a||b|+1}\big((\Phi^{\ast})_{[a,b]}|\biomega\big)=
-\big((\Phi^{\ast})_{b}|(id\otimes AD_{\ast
a})\biomega\big)$$ 
where we have defined $\big(\xi|(id\otimes AD_{\ast a})\big)
\comp\biomega\big)=(-1)^{|a||\xi|}
(id\otimes AD_{\ast a})(\xi|\biomega)$.

We may thus write:
$$\Lie_{(\Phi^{\ast})_{a}}\biomega
=-(id\otimes AD_{\ast a})\comp\biomega.\eqn\liead$$
Summarizing the previous results on the $\frak g$-valued graded
differential form $\biomega$, we have:
\item{1.}
$(\sum_{i}f^{i}(\Phi^{\ast})_{e_{i}}|\biomega)=
\sum_{i}f^{i}\otimes e_{i}$ and $|\biomega|=(1,0)$
\item{2.} $(\Phi^{\ast})_{a}\biomega=
(id\otimes AD_{\ast s^{\circ}_{\script
A}(a)})\comp\biomega$,
for each $a$ group-like or primitive with respect to $\delta_{e}$.

Conversely now, if $\biomega$ is a $\frak g$-valued graded
differential form with the previous two properties, then 
$\ker\biomega=\H$ is a regular distribution on $(Y,\B)$ such that
$\H\oplus\Ver(\pi_{\ast},\B)=\Der\B$. Consider $D\an\H(V)$,
$V=\pi_{\ast}^{-1}(U)$; if $a=\delta_{g}$ and $\biomega_{V}$ 
means the pull-back of $\biomega$ under the inclusion 
$(V,\B|_{V})\hookrightarrow(Y,\B)$, we take: 
$(D|\Phi_{g}^{\ast}\biomega_{V})=(\Phi_{g}^{\ast}\otimes
id)(\Phi_{g\ast}D|\biomega_{V})=
(id\otimes AD_{g^{-1}\ast})(D|\biomega_{V})=0$;
thus, $\Phi_{g\ast}D\an\H(V)$. Replacing $g$ by $g^{-1}$, this
gives $\Phi_{g}^{\ast}D\an\H(V)$. On the other hand, if $a\an\frak g$,
we find: $(D|\Lie_{(\Phi^{\ast})_{a}}\biomega_{V})=-(-1)^{|a||D|}
([(\Phi^{\ast})_{a},D]|\biomega_{V})=0$. This means that 
$\Lie_{(\Phi^{\ast})_{a}}D$ belongs also to $\H(V)$. In summary:
$(\Phi^{\ast})_{a}D\an\H(V)$. We have thus proved the proposition.

\math{Proposition.}{\connexiongequiv}{\sl The graded connection 
$\H$ of Definition {\connexiong} is described equivalently by a
$\frak g$-valued graded differential 1-form $\biomega\an
\Omega^{1}(Y,\B,\frak g)$ of total $\Z_{2}$ degree zero such that:
\item{1.}
$(\sum_{i}f^{i}(\Phi^{\ast})_{e_{i}}|\biomega)=
\sum_{i}f^{i}\otimes e_{i}$,
\item{2.} $(\Phi^{\ast})_{a}\biomega=
(id\otimes AD_{\ast s^{\circ}_{\script A}(a)})\comp\biomega,$
for each element $a$ group-like or primitive with respect to 
$\delta_{e}$.

We call $\biomega$ graded connection form.}

\vskip0.3cm

Graded principal bundles are by
definition locally isomorphic to products of open graded submanifolds
of the base space by the structure graded Lie group. So, it is quite
natural to ask how one can construct graded connections on the trivial
graded principal bundle of Example {\paradeigmaena}.

\math{Example.}{\paradeigmatria}{Let $(Y,\B)=(X,\C)\times(G,\A)$ be 
as in Example {\paradeigmaena}. The right action $\Phi$ of $(G,\A)$ 
on $(Y,\B)$ is such that $\Phi^{\ast}=id\otimes\da$. Consider now a 
$\frak g$-valued 1-form $\beta\an\Omega^{1}(X,\C,\frak g)$ and a 
basis $\{e_{k}\}$ of $\frak g$. One can write $\beta=\sum_{k}
\beta^{k}\otimes e_{k}$, $\beta^{k}\an\Omega^{1}(Y,\B)$ with
$|\beta|=(1,0)$. If $\xi\an\Der\C(X)$ and $\eta\an\Der\A(G)$, we 
define a $\frak g$-valued 1-form $\biomega\an\Omega^{1}(Y,\B,\frak g)$
as:
$$(\xi\otimes\ena_{\script A}+\ena_{\script C}\otimes\eta|\biomega)=
\sum_{k}(\xi|\beta^{k})\otimes((L^{\ast})_{e_{k}}|\theta)+\ena_{\script
C}\otimes(\eta|\theta).\eqn\syndesh$$
In the previous relation $L\Colon(G,\A)\times(G,\A)\rightarrow(G,\A)$ is
the left action of $(G,\A)$ on itself, $L^{\ast}=\da$ and
$\theta\an\Omega^{1}(G,\A,\frak g)$ the graded Maurer-Cartan form on
$(G,\A)$ defined as $((R^{\ast})_{a}|\theta)=\ena_{\script A}\otimes
a$, $\forall a\an\frak g$, where $R\Colon(G,\A)\times(G,\A)
\rightarrow(G,\A)$ is the right action of $(G,\A)$ on itself,
$R^{\ast}=\da$. Recall that the derivations $(L^{\ast})_{e_{k}}$ and
$(R^{\ast})_{a}$ are given by Theorem {\thewrima}. We shall check now
if the form $\biomega$ in $\syndesh$ has the properties of a
graded connection form. 

(1) If $a\an\frak g$, then $(\Phi^{\ast})_{a}=id
\otimes(R^{\ast})_{a}$ and $((\Phi^{\ast})_{a}|\biomega)=
\ena_{\script C}\otimes((R^{\ast})_{a}|\theta)=\ena_{\script B}\otimes
a$. Further, $|\biomega|=(1,0)$ because $|\beta|=(1,0)$ and
$|\theta|=(1,0)$.   

(2) Consider now an element $g\an G$; we shall calculate
$\Phi^{\ast}_{g}\biomega$. To this end, we use formula $\neostypos$,
as well as the fact that each graded Lie group is parallelizable
(Theorem {\paralleltheorem}), so it is sufficient to take
$\eta=(R^{\ast})_{a}$, $a\an\frak g$. Thus, if we set $D=\xi\otimes
\ena_{\script A}+\ena_{\script C}\otimes(R^{\ast})_{a}$ and
$a\tonos=AD_{g^{-1}\ast}(a)$, we find
$$\eqalign{(D|\Phi^{\ast}_{g}\biomega)&=
(\Phi^{\ast}_{g}\otimes id)(\xi\otimes\ena_{\script A}+\ena_{\script
C}\otimes(R^{\ast})_{a\tonoss}|\biomega)\cr
\hfill&=\sum_{k}(\xi|\beta^{k})\otimes(R^{\ast}_{g}\otimes
id)((L^{\ast})_{e_{k}}|\theta)+\ena_{\script B}\otimes
AD_{g^{-1}\ast}(a)\cr
\hfill&=\sum_{k}(\xi|\beta^{k})\otimes(id\otimes AD_{g^{-1}\ast})
((L^{\ast})_{e_{k}}|\theta)+\ena_{\script B}\otimes AD_{g^{-1}\ast}(a)\cr 
\hfill&=(\xi\otimes\ena_{\script A}+\ena_{\script
C}\otimes(R^{\ast})_{a}|(id\otimes
AD_{g^{-1}\ast})\comp\biomega).\cr}$$
Note that in the previous calculation we used that
$R^{\ast}_{g}\theta=(id\otimes AD_{g^{-1}\ast})\comp\theta$ and
$R_{g\ast}(L^{\ast})_{e_{k}}=(L^{\ast})_{e_{k}}$, which can be
verified straightforwardly. 

(3) Let finally $a\an\frak g$; we calculate the Lie derivative
$\Lie_{(\Phi^{\ast})_{a}}\biomega=(\Phi^{\ast})_{a}\biomega$. One can
use formula $\xrhsimo$, which is valid for $\biomega$ because
$((\Phi^{\ast})_{a}|\biomega)=\ena_{\script B}\otimes a$. If we set now
$D=\xi\otimes\ena_{\script A}+\ena_{\script C}\otimes(R^{\ast})_{b}$
and $x=|a||\beta^{k}|$, we have: 
$$\eqalign{(D|(\Phi^{\ast})_{a}\biomega)&=(-1)^{|a||\xi|}
((\Phi^{\ast})_{a}\otimes id)(\xi\otimes\ena_{\script
A}|\biomega)-(-1)^{|a||b|}(\ena_{\script
C}\otimes(R^{\ast})_{[a,b]}|\biomega)\cr
\hfill&=(-1)^{x}(\xi|\beta^{k})\otimes ((R^{\ast})_{a}\otimes id)
((L^{\ast})_{e_{k}}|\theta)-(-1)^{|a||b|}\ena_{\script B}\otimes[a,b]\cr
\hfill&=-(-1)^{x}(\xi|\beta^{k})\otimes(id\otimes AD_{\ast a})
((L^{\ast})_{e_{k}}|\theta)-(-1)^{|a||b|}\ena_{\script B}\otimes[a,b]\cr
\hfill&=-(\xi\otimes\ena_{\script A}|(id\otimes AD_{\ast
a})\comp\biomega)-(\ena_{\script C}\otimes(R^{\ast})_{b}|(id\otimes
AD_{\ast a})\comp\biomega)\cr}$$
which implies that $(\Phi^{\ast})_{a}\biomega=-(id\otimes AD_{\ast
a})\comp\biomega$. In the previous calculation, summation over the 
repeated index $k$ is understood. Note also that we used the property
$(R^{\ast})_{a}\theta=-(id\otimes AD_{\ast a})\comp\theta$ of the
graded Maurer-Cartan form. }

\vskip0.2cm

\math{Remark.}{\tmisi}{\it Let $s\Colon(X,\C)\rightarrow(Y,\B)$ be a
graded section for the previous example. Then we have a morphism 
$s^{\ast}\Colon\C(X)\otimes\A(G)\rightarrow\C(X)$ of graded commutative algebras  and therefore, a morphism $\sigma^{\ast}\Colon
\A(G)\rightarrow\C(X)$ given by $\sigma^{\ast}(f_{\script
A})=s^{\ast}(\ena_{\script C}\otimes f_{\script A})$. This defines 
a morphism of graded manifolds $\sigma\Colon(X,\C)\rightarrow(G,\A)$. 
Let now $\xi\an\Der\C(X)$ and $\eta\an\Der\A(G)$ be two
$\sigma$-related derivations; then the derivations
$\xi\tonos=\xi\otimes\ena_{\script A}+\ena_{\script C}\otimes\eta$ 
and $\xi$ are $s$-related and one can use Proposition {\pullback}
in order to calculate the pull-back $s^{\ast}\biomega\an\Omega^{1}
(X,\C,\frak g)$. The result is:
\vskip0.2cm
\centerline{$s^{\ast}\biomega=\sum_{k}\beta^{k}(\sigma^{\ast}\otimes
id)((L^{\ast})_{e_{k}}|\theta)+\sigma^{\ast}\theta,$}
\vskip0.1cm
\noindent as one easily finds.}

\vskip0.2cm

One can use Example {\paradeigmatria}  in order to construct graded
connections on general graded principal bundles. Indeed, let
$(Y,\B)$ be such a bundle with base space $(X,\C)$ and structure group
$(G,\A)$, $\{U_{i}\}_{i\an\Lambda}$ a locally finite open 
covering of $X$ and $\{f_{i}\}$ a graded partition of unity 
subordinate to $\{U_{i}\}$: $f_{i}\an\C(X)_{0}$, ${\rm supp}f_{i}
\subset U_{i}$ and $\sum_{i}f_{i}=\ena_{\script C}$, {\kost}. Let also
$V_{i}=\pi_{\ast}^{-1}(U_{i})$ and $\biomega_{i}$ be the graded
connection 1-forms that one can construct on 
$(V_{i},\B|_{V_{i}})\cong(U_{i},\C|_{U_{i}})\times(G,\A)$ as in the
previous example (see also Definition {\gpb}). If now $D\an\Der\B(Y)$,
then we define $(D|\biomega_{i})\an\B(V_{i})\otimes\frak g$ as
$(D|\biomega_{i})=(D_{i}|\biomega_{i})$, $D_{i}=D|_{V_{i}}$. Then we
have also $(D|\biomega_{i})\cdot(\pi^{\ast}f_{i}|_{V_{i}})\an\B(V_{i})
\otimes\frak g$ and so, there exists an element of $\B(Y)\otimes\frak
g$, which we denote by $(D|\biomega_{i})\pi^{\ast}f_{i}$, such that 
$[(D|\biomega_{i})\pi^{\ast}f_{i}]|_{V_{i}}=(D|\biomega_{i})
(\pi^{\ast}f_{i}|_{V_{i}})$ and whose support is a subset of $V_{i}$.

If now $\{\tilde{U}_{k}\}$ is another open covering of $X$, then 
setting $W_{k}=\pi_{\ast}^{-1}(\tilde{U}_{k})$ we have an open
covering of $Y$ and for $k$ fixed, $W_{k}\cap V_{i}$
is non-empty only for finitely many of the $V_{i}$. Taking the
restrictions $[(D|\biomega_{i})\pi^{\ast}f_{i}]|_{W_{k}}\an\B(W_{k})
\otimes\frak g$, only finite many terms will be non-zero as $i$ runs
over $\Lambda$. Thus, the sum $\sum_{i}[(D|\biomega_{i})
\pi^{\ast}f_{i}]|_{W_{k}}$ is finite and well-defined as an element of
$\B(W_{k})\otimes\frak g$. Furthermore, the restrictions of such
elements to the intersections $W_{k}\cap W_{\ell}$ coincide and this
means, by the sheaf properties of $\B$, that there exists a unique
element of $\B(Y)\otimes\frak g$, whose restriction to $W_{k}$ gives
the previous element of $\B(W_{k})\otimes\frak g$. We denote this
unique element by $(D|\biomega)$ and by its linearity on $D$ it
determines a $\frak g$-valued graded 1-form $\biomega$. The form
$\biomega$ is a graded connection form; this results without 
difficulty from the properties $\Phi_{g}^{\ast}\pi^{\ast}f_{i}=
\pi^{\ast}f_{i}$ and $(\Phi^{\ast})_{a}\pi^{\ast}f_{i}=0$, $g\an G$, 
$a\an\frak g$, as well as from the fact that for each $i\an\Lambda$, 
$\biomega_{i}$ is a graded connection. We have thus proved the following:

\math{Existence Theorem for graded connections.}{\yparxh}{\sl On each
graded principal bundle there exists an infinity of graded
connections.}

\vskip0.2cm

\chapter{Graded curvature}

In ordinary differential geometry, one can define in a canonical way,
for each connection, a Lie superalgebra-valued 2-form, the curvature
form. In this section, we will define the curvature in the graded
setting, using the notion of graded connection, previously developed. 

Let then $\biomega\an\Omega^{1}(Y,\B,\frak g)$ be a graded connexion form on
the graded principal bundle $(Y,\B)$. Fix a basis of the Lie 
$\{e_{k}\}$ of the Lie superalgebra $\frak g$; then for each
derivation $\xi\an\Der\B(Y)$ one can find $\xi^{H}\an\H(Y)$ and 
$f^{k}\an\B(Y)$ such that $\xi=\xi^{H}+\sum_{i}f^{k}
(\Phi^{\ast})_{e_{k}}$. We have thus a canonical projection 
$\Der\B(Y)\rightarrow\H(Y)$, $\xi\mapsto\xi^{H}$;
we call $\xi^{H}$ horizontal part of $\xi$. Clearly, this 
mechanism of taking the horizontal part of a derivation can be 
applied in the same way if instead of $Y$ we put an open $V\subset Y$. 

Consider now a $\frak g$-valued graded differential form 
$\phi\an\Omega^{r}(Y,\B,\frak g)$ and let $\phi^{H}$ be 
defined as $(\xi_{1},\ldots,\xi_{r}|\phi^{H})=
(\xi_{1}^{H},\ldots,\xi_{r}^{H}|\phi)$,
$\xi_{i}\an\Der\B(Y)$, $i=1,\ldots,r$.

\math{Definition.}{\covariante}{\sl The covariant exterior derivative
of $\phi\an\Omega^{r}(Y,\B,\frak g)$ is the $\frak g$-valued graded 
differential form $D^{\biomega}\phi\an\Omega^{r+1}(Y,\B,\frak g)$ 
defined as $D^{\biomega}\phi=(d\phi)^{H}$. 
The curvature of the graded connection $\biomega\an
\Omega^{1}(Y,\B,\frak g)$ is the covariant exterior 
derivative $D^{\biomega}\biomega$; we use the notation 
$F^{\biomega}=D^{\biomega}\biomega$.}

\vskip0.3cm

We will next study in more detail the properties of $F^{\biomega}$. 
As a general observation, we may say that $F^{\biomega}$ has,
formally, properties analogous to those of the ordinary curvature.

\math{Theorem.}{\structureth}{\sl The graded curvature
$F^{\biomega}$ is given by the relation 
$$F^{\biomega}=d\biomega+{{1}\over{2}}[\biomega,\biomega].
\eqn\structureequ$$
This is the graded structure equation.}

\undertext{\it Proof}. It is sufficient to prove that for
each $\xi_{1}$, $\xi_{2}\an\Der\B(Y)$, the following is true:
$$(\xi_{1}^{H},\xi_{2}^{H}|d\biomega)=(\xi_{1},\xi_{2}|d\biomega)+
{{1}\over{2}}(\xi_{1},\xi_{2}|[\biomega,\biomega]).\eqn\structureequdyo$$
As we have seen $|\biomega|=(1,0)$; however, $\biomega$ is not a
homogeneous element with respect to the $(\Z\oplus\Z_{2}\oplus
\Z_{2})$-grading of $\OYBG$. More precisely, $\biomega=\biomega_{0}+
\biomega_{1}$, where $\deg(\biomega_{0})=(1,0,0)$, 
$\deg(\biomega_{1})=(1,1,1)$ and if $\{e_{i},e_{j}\}$ is a 
basis of $\frak g$ with $|e_{i}|=0$, $|e_{j}|=1$, we may write:
$\biomega=\sum_{i}\biomega^{i}\otimes e_{i}+\sum_{j}\biomega^{j}
\otimes e_{j}$, with $\biomega^{i}\an\Omega^{1}(Y,\B)_{0}$, $\biomega^{j}
\an\Omega^{1}(Y,\B)_{1}$. In particular, $\biomega^{i}$, $\biomega^{j}$
vanish on horizontal derivations and if we decompose $a\an\frak g$ as
$a=\sum_{i}a^{i}e_{i}+\sum_{j}a^{j}e_{j}$, we find immediately
$((\Phi^{\ast})_{a}|\biomega^{i})=a^{i}\ena_{\script B}$, 
$((\Phi^{\ast})_{a}|\biomega^{j})=a^{j}\ena_{\script B}$.
Let us now calculate the term ${{1}\over{2}}(\xi_{1},\xi_{2}|
[\biomega,\biomega])$. 
$$\eqalign{(\xi_{1},\xi_{2}|[\biomega,\biomega])&=
\big(\xi_{1},\xi_{2}\big|\big[{\textstyle\sum_{i}}\biomega^{i}\otimes
e_{i}+{\textstyle\sum_{j}}\biomega^{j}\otimes e_{j},
{\textstyle\sum_{p}}\biomega^{p}\otimes
e_{p}+{\textstyle\sum_{q}}\biomega^{q}\otimes e_{q}\big]\big)\cr
\hfill&={\textstyle\sum_{i,p}}
\Big((\xi_{1}|\biomega^{i})(\xi_{2}|\biomega^{p})+
(-1)^{1+|\xi_{1}||\xi_{2}|}(\xi_{2}|\biomega^{i})
(\xi_{1}|\biomega^{p})\Big)
\otimes[e_{i},e_{p}]\cr
\hfill&+{\textstyle\sum_{i,q}}\Big((\xi_{1}|\biomega^{i})
(\xi_{2}|\biomega^{q})+(-1)^{1+|\xi_{1}||\xi_{2}|}(\xi_{2}|\biomega^{i})
(\xi_{1}|\biomega^{q})\Big)\otimes[e_{i},e_{q}]\cr
\hfill&+{\textstyle\sum_{j,p}}
\Big((-1)^{|\xi_{2}|}(\xi_{1}|\biomega^{j})
(\xi_{2}|\biomega^{p})+(-1)^{x}
(\xi_{2}|\biomega^{j})(\xi_{1}|\biomega^{p})\Big)\otimes[e_{j},e_{p}]\cr
\hfill&-{\textstyle\sum_{j,q}}
\Big((-1)^{|\xi_{2}|}(\xi_{1}|\biomega^{j})
(\xi_{2}|\biomega^{q})+(-1)^{x}
(\xi_{2}|\biomega^{j})(\xi_{1}|\biomega^{q})\Big)\otimes[e_{j},e_{q}].\cr}
\eqn\makrinari$$
In the previous calculation, the indices $i,p$ label the even elements
while $j,q$ the odd ones. We also used the fact that if
$\beta_{1},\beta_{2}\an\Omega^{1}(Y,\B)$, then 
$(\xi_{1},\xi_{2}|\beta_{1}\beta_{2})=
(-1)^{|\xi_{2}||\beta_{1}|}(\xi_{1}|\beta_{1})(\xi_{2}|\beta_{2})+
(-1)^{1+|\xi_{1}||\xi_{2}|+|\xi_{1}||\beta_{1}|}(\xi_{2}|\beta_{1})
(\xi_{1}|\beta_{2})$, see relation 4.1.9 of {\kost}, and we have set 
$x=1+|\xi_{1}||\xi_{2}|+|\xi_{1}|$. On the other hand, the term 
$(\xi_{1},\xi_{2}|d\biomega)$ is calculated via 
$(\xi_{1},\xi_{2}|d\biomega)=(\xi_{1}\otimes id)(\xi_{2}|\biomega)-
(-1)^{|\xi_{1}||\xi_{2}|}(\xi_{2}\otimes id)(\xi_{1}|\biomega)-
([\xi_{1},\xi_{2}]|\biomega)$, which is an immediate generalization of
4.3.10 of {\kost}. We distinguish now the following cases:

\item{1.} $\xi_{1}$, $\xi_{2}$: horizontal
$\Rightarrow\xi_{1}=\xi_{1}^{H},\xi_{2}=\xi_{2}^{H}$.
The graded structure equation holds, since
${1\over 2}(\xi_{1},\xi_{2}|[\biomega,\biomega])=0$ and 
$(\xi_{i}|\biomega)=0$.
\item{2.} $\xi_{1}$: horizontal, $\xi_{2}=(\Phi^{\ast})_{a}$, 
$a\an\frak g$. The left-hand side of $\structureequdyo$ is zero
because $\xi_{2}^{H}=0$. The right-hand side of the same equation
reads: $(\xi_{1},(\Phi^{\ast})_{a}|d\biomega)=
(\xi_{1}\otimes id)((\Phi^{\ast})_{a}|\biomega)-([\xi_{1},
(\Phi^{\ast})_{a}]|\biomega)=0$, because by the definition of the
graded connexion, if $\xi$ is horizontal, $[\xi,(\Phi^{\ast})_{a}]$ is
horizontal too. Furthermore, it is clear that in this case,
$\makrinari$ gives $(\xi_{1},\xi_{2}|[\biomega,\biomega])=0$.

\item{3.} $\xi_{1}=(\Phi^{\ast})_{a}$, $\xi_{2}=(\Phi^{\ast})_{b}$,
$a,b\an\frak g$. Clearly, the left-hand side of $\structureequdyo$ is
zero. Examining now the cases $|a|=|b|=0$, $|a|=|b|=1$ and
$|a|=1,|b|=0$, we find always that the right-hand side is also
zero.\qed

\vskip0.3cm

By the graded structure equation, it is clear that $\biomega$ and
$F^{\biomega}$ have the same $\Z_{2}$ total degree:
$|F^{\biomega}|=(2,0)$; therefore,
$F^{\biomega}=F^{\biomega}_{0}+F^{\biomega}_{1}$ with 
$F^{\biomega}_{0}\an\Omega^{2}(Y,\B)_{0}\otimes\frak g_{0}$, 
$F^{\biomega}_{1}\an\Omega^{2}(Y,\B)_{1}\otimes\frak g_{1}$. 

Bianchi's identity is a well-known property satisfied by the curvature
in differential geometry. Using the generalization of the covariant
derivative in the graded setting and the previous theorem, we may 
establish an analogous property in the context of graded manifolds. 

\math{Proposition. (Bianchi's Identity)}{\bianchi}
{$D^{\biomega}F^{\biomega}=0$.}  

\undertext{\it Proof}. Let us first calculate the differential 
$dF^{\biomega}$. Using the graded structure equation and the fact that
$\big[[\biomega,\biomega],\biomega\big]=0$ (Jacobi identity), we find easily:
$$dF^{\biomega}={1\over 2}d[\biomega,\biomega]
={1\over 2}([d\biomega,\biomega]-[\biomega,d\biomega])
=[d\biomega,\biomega]=[F^{\biomega},\biomega].$$
Thus, $dF^{\biomega}=[F^{\biomega},\biomega]$ and
$(\xi_{1},\xi_{2},\xi_{3}|D^{\biomega}F^{\biomega})=
(\xi_{1}^{H},\xi_{2}^{H},\xi_{3}^{H}|[F^{\biomega},\biomega])=0$,
because $\biomega$ vanishes on horizontal derivations.\qed

\vskip0.3cm

We will show now that the graded curvature $F^{\biomega}$ satisfies 
the second property of the connexion $\biomega$ described in
Proposition {\connexiongequiv}. To this end, consider first 
$a=\delta_{g}$, $g\an G$. Then:
$$\eqalign{\Phi^{\ast}_{g}F^{\biomega}&=\Phi^{\ast}_{g}
(d\biomega+{1\over 2}[\biomega,\biomega])
=d\Phi^{\ast}_{g}\biomega+{1\over
2}\Phi^{\ast}_{g}[\biomega,\biomega]\cr
\hfill&=(id\otimes AD_{g^{-1}\ast})d\biomega+{1\over 2}
(id\otimes AD_{g^{-1}\ast})[\biomega,\biomega]
=(id\otimes AD_{g^{-1}\ast})F^{\biomega}.\cr}$$ 
Suppose now that $a\an\frak g$ homogeneous; thanks to Proposition
{\idiothtalie}, direct calculation gives:
$$\eqalign{(\Phi^{\ast})_{a}F^{\biomega}&=
\Lie_{(\Phi^{\ast})_{a}}F^{\biomega}
=\Lie_{(\Phi^{\ast})_{a}}d\biomega+{1\over
2}\Lie_{(\Phi^{\ast})_{a}}[\biomega,\biomega]\cr
\hfill&=d\Lie_{(\Phi^{\ast})_{a}}\biomega+{1\over 2}
\big([\Lie_{(\Phi^{\ast})_{a}}\biomega,\biomega]+
[\biomega,\Lie_{(\Phi^{\ast})_{a}}\biomega]\big)\cr
\hfill&=-d(id\otimes AD_{\ast a})\biomega+
{1\over 2}\big([-(id\otimes AD_{\ast a})\biomega,\biomega]+
[\biomega,-(id\otimes AD_{\ast a})\biomega]\big)\cr
\hfill&=-\big((id\otimes AD_{\ast a})d\biomega+{1\over 2}
(id\otimes AD_{\ast a})[\biomega,\biomega]\big)
=-(id\otimes AD_{\ast a})F^{\biomega},\cr}$$
because $|\biomega|=(1,0)$, $|\Lie_{(\Phi^{\ast})_{a}}|=(0,|a|)$. We
have thus proved the following:

\math{Property.}{\toidioisxuei}{\sl If $a\an\A(G)^{\circ}$ is a
group-like or primitive element with respect to $\delta_{e}$, we have 
$(\Phi^{\ast})_{a}F^{\biomega}=(id\otimes AD_{\ast s^{\circ}_{\script
A}(a)})\comp F^{\biomega}$.}

\vskip0.3cm

We finally prove that the graded curvature provides a criterion for
checking whether the horizontal distribution is involutive or not.

\math{Theorem.}{\involutiveH}{\sl The graded horizontal distribution 
is involutive if and only if the graded curvature $F^{\biomega}$ of 
the connexion $\biomega$ is zero: 
$[\H,\H]\subset\H\Leftrightarrow F^{\biomega}=0$.}

\undertext{\it Proof}. It is sufficient to check the involutivity
of $\H(Y)$. If $\xi,\eta\an\H(Y)$ and $\H$ is involutive, then:
$([\xi,\eta]|\biomega)=0\Rightarrow(\xi,\eta|d\biomega)+{1\over
2}(\xi,\eta|[\biomega,\biomega])=0\Rightarrow(\xi,\eta|F^{\biomega})=0$.
But $F^{\biomega}$ vanishes identically on vertical derivations by its
definition (relation $\structureequdyo$); thus $F^{\biomega}=0$. 

Conversely, suppose that $F^{\biomega}=0$. Then for each 
$\xi,\eta\an\H(Y)$, we find: $(\xi,\eta|F^{\biomega})=0\Rightarrow 
(\xi,\eta|d\biomega)+{1\over 2}(\xi,\eta|[\biomega,\biomega])=0
\Rightarrow(\xi,\eta|d\biomega)=0\Rightarrow-([\xi,\eta]|\biomega)=0$,
which implies that the derivation $[\xi,\eta]$ is also horizontal, 
$[\xi,\eta]\an\H(Y)$.\qed

\vskip0.3cm

\chapter{Concluding remarks}

Consider a graded principal bundle $(Y,\B)$ equipped with a connection
form $\biomega\an\Omega^{1}(Y,\B,\frak g)$. We know from the general
theory of graded differential forms, {\kost}, that there always exists
an algebra morphism $\kappa\Colon\Omega(Y,\B)\rightarrow\Omega(Y)$
defined as follows: if $i\Colon(Y,C^{\infty})\rightarrow(Y,\B)$ is 
the morphism of
graded manifolds determined by $\B(Y)\na f\mapsto\tilde{f}\an
C^{\infty}(Y)$ (see exact sequence $\exactsequence$), then $i^{\ast}$
is just $\kappa$. On the other hand, the decomposition $\frak g=\frak
g_{0}\oplus\frak g_{1}$ induces a canonical projection
$\epi_{0}\Colon\frak g\rightarrow\frak g_{0}$. So we have a linear map
$\kappa_{0}=\kappa\otimes\epi_{0}\Colon\Omega(Y,\B,\frak
g)\rightarrow\Omega(Y)\otimes\frak g_{0}$. Explicitly, if
$\alpha=\sum_{i}\alpha^{i}\otimes e_{i}\an\Omega(Y,\B,\frak g)$, then 
$\kappa_{0}(\alpha)=\sum_{i}\kappa\big((\alpha^{i})_{0}\big)
\otimes(e_{i})_{0}$, where $(\alpha^{i})_{0}\an\Omega(Y,\B)_{0}$ and 
$(e_{i})_{0}\an\frak g_{0}$ are the even elements in the development
of $\alpha$; moreover, we easily realize that $\kappa_{0}$ is not a 
$(\Z\oplus\Z_{2})$-graded Lie algebra morphism. 

As we have seen in the proof of Theorem {\structureth}, $\biomega$ can
be decomposed as $\biomega=\biomega_{0}+\biomega_{1}$, where 
$\biomega_{0}\an\Omega^{1}(Y,\B)_{0}\otimes\frak g_{0}$ and
$\biomega_{1}\an\Omega^{1}(Y,\B)_{1}\otimes\frak g_{1}$. Then clearly,
$\kappa_{0}(\biomega)=\kappa_{0}(\biomega_{0})=\sum_{i}\kappa
(\biomega^{i})\otimes e_{i}$, where $i$ labels the even elements.
Using the fact that the derivations induced on $(Y,\B)$ and
$(Y,C^{\infty})$ by the right actions of $(G,\A)$ and $(G,C^{\infty})$
respectively (according to Theorem {\thewrima}) are $i$-related, as
well as the defining properties of a graded connection form, one can
prove that $\kappa_{0}(\biomega)$ is a connection form on the ordinary
principal bundle $(Y,C^{\infty})$. Furthermore, the curvatures of
$\biomega$ and $\kappa_{0}(\biomega)$ are related through $\kappa_{0}
(F^{\biomega})=F^{\kappa_{0}(\biomega)}$. 

The previous observations suggest that the connection theory on graded
principal bundles is the suitable framework for the mathematical
formulation of super-gauge field theories. The fact that the graded
connection $\biomega$ splits always as $\biomega=\biomega_{0}+
\biomega_{1}$ with $\kappa_{0}(\biomega_{0})$ being a usual connection
form, incorporates automatically the idea of supersymmetric partners:
$\biomega_{0}$ corresponds to the gauge potential of an ordinary
Yang-Mills theory, while $\biomega_{1}$ corresponds to its
supersymmetric partner. A very interesting feature of this approach is
that there is no graded connection form $\biomega$ for which one of
the terms $\biomega_{0}$ or $\biomega_{1}$ is zero. In physics
terminology, all gauge potentials have super-partners. The same is
also true for the curvature $F^{\biomega}$, the super-gauge field, 
since $F^{\biomega}=F^{\biomega}_{0}+F^{\biomega}_{1}$,
$F^{\biomega}_{0}\an\Omega^{2}(Y,\B)_{0}\otimes\frak g_{0}$,  
$F^{\biomega}_{1}\an\Omega^{2}(Y,\B)_{1}\otimes\frak g_{1}$.
Thus, our approach sets the gauge potentials (resp. fields) and
its supersymmetric partners on the same footing: they both ``live"
in a Lie superalgebra $\frak g$ as components of the same connection
form (resp. curvature form). This is an essential difference between 
our approach and the standard treatment of this problem by means of 
DeWitt's or Roger's supermanifolds, (see {\grasso} and references
therein), where the connections take values in the even part of the
$\Z_{2}$-graded Lie module corresponding to a Lie supergroup.

\vskip0.5cm

\noindent {\bf Acknowledgements.} 

\vskip0.1cm

I would like to thank Professor R. Coquereaux 
for his critical reading of the manuscript and for many stimulating
discussions.

\vskip1cm
   \ifreferenceopen \Closeout\referencewrite \referenceopenfalse \fi
   \line{\bf\hskip0pt\hfil References\hfil}\vskip\headskip
   \vskip0.3cm
   \input referenc.txa
\end